\documentclass[11pt]{article}
 \usepackage{benstyle}
 
\newtheorem{exampleresult}{Informal Result}

\title{Learning smooth functions in high dimensions: from sparse polynomials to deep neural networks}

\author{Ben Adcock \thanks{Department of Mathematics, Simon Fraser University, 8888 University Drive, Burnaby BC, Canada, V5A 1S6. \emph{E-mail address:} {\tt ben\_adcock@sfu.ca} } \and Simone Brugiapaglia\thanks{Department of Mathematics and Statistics, Concordia University, J.W. McConnell Building, 1400 De Maisonneuve Blvd. W., Montréal, QC, Canada, H3G 1M8.  \emph{E-mail address:} {\tt simone.brugiapaglia@concordia.ca} } \and Nick Dexter\thanks{Department of Scientific Computing, Florida State University,
400 Dirac Science Library, Tallahassee, Florida, USA, 32306-4120.  \emph{E-mail address:}  {\tt nick.dexter@fsu.edu} } \and Sebastian Moraga \thanks{Department of Mathematics, Simon Fraser University, 8888 University Drive, Burnaby BC, Canada, V5A 1S6. 
Corresponding author: {\tt smoragas@sfu.ca}  }}

\begin{document}

\maketitle

\begin{abstract}
Learning approximations to smooth target functions of many variables from finite sets of pointwise samples is an important task in scientific computing and its many applications in computational science and engineering. Despite well over half a century of research on high-dimensional approximation, this remains a challenging problem. Yet, significant advances have been made in the last decade towards efficient methods for doing this, commencing with so-called \textit{sparse polynomial approximation} methods and continuing most recently with methods based on \textit{Deep Neural Networks (DNNs)}. In tandem, there have been substantial advances in the relevant approximation theory and analysis of these techniques. In this work, we survey this recent progress. We describe the contemporary motivations for this problem, which stem from parametric models and computational uncertainty quantification; the relevant function classes, namely, classes of infinite-dimensional, Banach-valued, holomorphic functions; fundamental limits of learnability from finite data for these classes; and finally, sparse polynomial and DNN methods for efficiently learning such functions from finite data. For the latter, there is currently a significant gap between the approximation theory of DNNs and the practical performance of deep learning. Aiming to narrow this gap, we develop the topic of \textit{practical existence theory}, which asserts the existence of dimension-independent DNN architectures and training strategies that achieve provably near-optimal generalization errors in terms of the amount of training data.
\end{abstract}


\section{Introduction}\label{s:introduction}

This work reviews the problem of learning infinite-dimensional functions from finite data. In this section, with begin by describing motivations for this problem and its challenges (\S \ref{ss:motivations}). We then present an overview of the remainder of the work (\S \ref{s:overview}). Finally, we discuss relevant literature (\S\ref{ss:literature}).

\subsection{Motivations and challenges}\label{ss:motivations}

Approximating functions of many variables is a classical topic in approximation theory and numerical analysis, which has been studied intensively for more than seventy years. The contemporary resurgence in this problem stems from applications to parametric models, parametric Differential Equations (DEs) and computational Uncertainty Quantification (UQ). In such models, the target function represents some quantity-of-interest of a given physical system and its variables the parameters of the system. Complex physical models involve many parameters, which naturally  leads to functions of many variables. It is also increasingly common to consider countably infinite parameterizations, a standard example being the use of a Karhunen--Lo\`eve expansion  to represent a random field in, for example, parametric models of porous  media (see, e.g., \cite[\S 2.1]{le-maitre2010spectral} or \cite[\S 11.1]{sullivan2015introduction}). In this case, the target function depends on countably many variables.

The first challenge in approximating such functions, therefore, is that the dimension is high or countably infinite. The second is that the function is usually expensive to evaluate. Typically, each evaluation involves either an expensive numerical simulation of the DE model or a costly physical experiment. Therefore, one strives to learn an approximation from a limited set of samples. This approximation -- or \textit{surrogate model} as it is often termed (see \cite[Chpt.\ 13]{smith2013uncertainty} or \cite[Chpt.\ 13]{sullivan2015introduction}) -- can then be used in place of the true function for key downstream tasks, such as parameter optimization, inverse parametric problems or uncertainty quantification. Nonetheless, since the target function is expensive to evaluate, the amount of data available to learn this approximation is usually highly limited. Put another way, the learning problem is highly \textit{data-starved}. 

A final challenge is that the target function is often not scalar-valued. In the case of parametric DEs, the output of the target function is the solution of the DE with the given input parameters. Hence this function takes values in infinite-dimensional  Banach or Hilbert space. This adds complications, both theoretically and practically, stemming from the need to discretize this output space in order to perform computations.

\subsection{Overview}\label{s:overview}

We now give a detailed overview of this article.

\subsubsection*{Problem statement and notation (\S \ref{s:prob-statement})}

We consider an unknown target function depending on countably many variables that range between finite maximum and minimum values and taking values in a Banach space $(\cV,\nm{\cdot}_{\cV})$. After normalizing, we may assume that
\bes{
f : \cU \rightarrow \cV,\quad \text{where }\cU = [-1,1]^{\bbN}.
}
We equip $\cU$ with the uniform probability measure $\varrho$ and consider training data
\be{
\label{training-data-intro}
\{ (\bm{y}_i,f(\bm{y}_i) + e_i ) \}^{m}_{i=1} \subseteq \cU \times \cV,
}
where $\bm{y}_1,\ldots,\bm{y}_m \sim_{\mathrm{i.i.d.}} \varrho$ and $e_i \in \cV$ represents measurement error. The objective is to learn an approximation $\hat{f}$ to $f$ from the data \ef{training-data}.

\subsubsection*{Holomorphic functions of infinitely many variables (\S \ref{s:holo-fns})}

An important property of many parametric model problems is that the parametric map $f$ is a smooth function of its variables. There is a large body of literature, which we review further in \S \ref{s:holo-fns}, that establishes that solution maps for various parametric DEs are \textit{holomorphic} functions -- in other words, they admit holomorphic extensions to certain complex regions $\cU \subset \cO \subseteq \bbC^{\bbN}$.

In one dimension, it is common to consider complex regions $[-1,1] \subset \cO \subset \bbC$ defined by \textit{Bernstein ellipses} (see, e.g., \cite{trefethen2013approximation}). The convergence of the $s$-term polynomial expansion is then \textit{exponential} in $s$, with the rate being characterized by the largest Bernstein ellipse to which the function can be extended (see, e.g., \cite[Chpt.\ 8]{trefethen2013approximation}). For functions of $d \geq 1$ variables, it is natural to look at regions $[-1,1]^d \subset \cO \subset \bbC^d$ defined by \textit{Bernstein polyellipses}, i.e., tensor-products of one-dimensional Bernstein ellipses. One can then show exponential convergence in $s^{1/d}$ of a suitable $s$-term polynomial expansion, with rate depending once more on the largest polyellipse to which $f$ admits an extension (see, e.g., \cite[\S 3.5-3.6]{adcock2022sparse}).

\pbk
\textit{$(\bm{b},\varepsilon)$-holomorphic functions.} The situation changes in infinite dimensions. Holomorphy in an arbitrary Bernstein polyellipse is no longer sufficient to guarantee convergence of a (in fact, any) $s$-term polynomial expansion. At the very least, one needs \textit{anisotropy}, namely, increasing regularity with the variables $y_i$ as $i \rightarrow \infty$. Moreover, for parametric DE applications, one needs to consider not just a single Bernstein ellipse, but complex regions $\cR(\bm{b},\varepsilon) \subset \bbC^{\bbN}$ defined by certain \textit{unions} of Bernstein polyellipses. Here $\bm{b} = (b_{j})_{j \in \bbN} \in [0,\infty)^{\bbN}$ with $\bm{b} \in \ell^1(\bbN)$ is a sequence and $\varepsilon > 0$ is a scalar that parametrize the corresponding region. In Definition \ref{d:b-eps-holomorphy} we formalize the corresponding class of so-called \textit{$(\bm{b},\varepsilon)$-holomorphic functions} (see, e.g., \cite{chkifa2015breaking,schwab2019deep}). This class was first introduced in the context of parametric DEs \cite{cohen2015approximation}, and has since become a standard setting in which to consider the approximation of holomorphic functions in infinite dimensions \cite[Chpt.\ 3]{adcock2022sparse}. Throughout the majority of this work, we set $\varepsilon = 1$; a property which can always be guaranteed by rescaling $\bm{b}$. For convenience, we also define the set
\be{
\label{Hb-intro}
\cH(\bm{b}) =  \big\{ f : \cU \rightarrow \cV\text{ $(\bm{b},1)$-holomorphic} : \sup_{\bm{z} \in \cR(\bm{b},1)} \nm{f({\bm{z}})}_{\cV}  \leq 1 \big\}.
}
This is the main class of functions considered in the remainder of this work. 
In \S \ref{s:holo-fns} we explain the relevance of this class to parametric DEs in more detail.

\pbk
\textit{Unknown anisotropy.} We next explain how the sequence $\bm{b}$ controls the \textit{anisotropy} of $f \in \cH(\bm{b})$. Larger $b_i$ implies that $f$ is less smooth with respect to its $i$th variable $y_i$, while smaller $b_i$ implies more smoothness. In some specific situations, one may have knowledge of a suitable $\bm{b}$ -- information which could then be used when designing a learning algorithm. However, in practice, it is more often the case that $\bm{b}$ is unknown. The focus of this work is the more realistic \textit{unknown anisotropy} setting. In short, the learning method must be independent of $\bm{b}$, with the assumption $f \in \cH(\bm{b})$ being used only  to derive bounds for the resulting generalization error.

\subsubsection*{Orthogonal polynomials and best $s$-term polynomial approximation (\S \ref{s:sparse-poly})}

As previously observed, holomorphy is intimately related to polynomial approximation. In \S \ref{s:sparse-poly} we introduce multivariate orthogonal polynomials, orthogonal polynomial expansions of infinite-dimensional, Banach-valued functions and the concept of \textit{best $s$-term polynomial approximation}. As we argue, best $s$-term polynomial approximation is an important theoretical benchmark against which to compare methods for learning holomorphic functions from data.

A signature result we recap in this section is the following.

\begin{exampleresult}[Algebraic convergence of the best $s$-term approximation, Theorem \ref{t:best-s-term}]
\textit{The best $s$-term polynomial approximation of $f \in \cH(\bm{b})$ converges algebraically fast in $s$. Specifically, if $\bm{b} \in \ell^p(\bbN)$ for some $0 < p < 1$, then the $L^2_{\varrho}$-norm error decreases like $\ord{s^{1/2-1/p}}$ as $s \rightarrow \infty$.}
\end{exampleresult}

This result implies that best $s$-term polynomial approximation in $\cH(\bm{b})$ is free from the \textit{curse of dimensionality}, whereas the aforementioned exponential rates witnessed in $d < \infty$ dimensions (which are exponential in $s^{1/d}$) quickly succumb to the curse of dimensionality. Although the focus in this work is infinite-dimensional functions, in \S \ref{s:sparse-poly} we demonstrate that the algebraic rates also typically describe the true convergence behaviour seen for functions of finitely many variables, whenever $d$ is large enough (e.g., $d \geq 16$).

\subsubsection*{Limits of learnability from data (\S \ref{s:learn-limits})}

We next turn our attention to learning holomorphic functions from the data \ef{training-data-intro}. In \S \ref{s:learn-limits}, we study lower bounds using concepts from \textit{information-based complexity}, in particular, \textit{adaptive $m$-widths} (Definition \ref{def:m-width}). In order to consider the unknown anisotropy setting, we introduce the function classes
\bes{
\cH(p) =  \bigcup \left \{ \cH(\bm{b}) : \bm{b} \in \ell^p(\bbN),\ \bm{b} \in [0,\infty)^{\bbN},\   \nm{\bm{b}}_p \leq 1 \right \} ,\quad 0 < p < 1,
}
where $\cH(\bm{b})$ is as in \ef{Hb-intro}, and 
\bes{
\cH(p,{\mathsf{M}})  = \bigcup  \left \{ \cH(\bm{b}) : \bm{b} \in \ell^p_{\mathsf{M}}(\bbN),\ \bm{b} \in [0,\infty)^{\bbN},\  \nm{\bm{b}}_{p,\mathsf{M}} \leq 1 \right \} ,\quad 0 < p < 1.
}
Here $\ell^p_{\mathsf{M}}(\bbN)$ is the \textit{monotone} $\ell^p$-space (see \S \ref{s:prob-statement}). We then have the following.

\begin{exampleresult}
[Limits of learnability, Theorem \ref{t:lower-bounds}]
\textit{It is impossible to learn functions from $\cH(p)$. Specifically, there does not exist a method for learning functions from $m$ (adaptive) linear samples for which the $L^2_{\varrho}$-norm error decreases as $m \rightarrow \infty$ uniformly for functions in $\cH(p)$. Moreover, when restricting to the space $\cH(p,\mathsf{M})$, or even $\cH(\bm{b})$, the error cannot decrease faster that $c \cdot m^{1/2-1/p}$ even if the method is allowed to depend on the anisotropy parameter $\bm{b} \in \ell^{p}_{\mathsf{M}}(\bbN)$.}
\end{exampleresult}

This theorem illustrates a fundamental gap between approximation theory and learning from data. The best $s$-term approximation converges like $s^{1/2-1/p}$ for any $f \in \cH(p)$, yet no method can learn such functions from data. Note that this result holds not just for i.i.d.\ pointwise samples, as in \ef{training-data-intro}, but \textit{any} linear measurements, which may also be \textit{adaptive}.

This gap is narrowed by restricting to the class $\cH(p,\mathsf{M})$, in which the variables are ``on average'' ordered in terms of importance. Notably, this lower bound also implies that there is no benefit to knowing $\bm{b}$, at least in terms of sample complexity. The same lower bound $c \cdot m^{1/2-1/p}$ holds, even if the method is allowed to depend on $\bm{b}$.

With this in mind, the remainder of this work is devoted to describing learning methods that achieve (close to) the rate $m^{1/2-1/p}$ uniformly for functions in $\cH(p,\mathsf{M})$ for \textit{any} $0 < p < 1$. Specifically, we focus on sparse polynomial and Deep Neural Network (DNN) methods.

\subsubsection*{Learning sparse polynomial approximations from data (\S\ref{s:sparse-poly-learn})}

Inspired by the benchmark results on best $s$-term approximation, we first consider methods that learn polynomial approximations. This approach is heavily based on techniques from \textit{compressed sensing} \cite{foucart2013mathematical,vidyasagar2019introduction}. Superficially, this seems straightforward. Fast decay of the best $s$-term polynomial approximation means that the polynomial coefficients are approximately sparse. Thus, a first approach would be to formulate an $\ell^1$-minimization problem for the coefficients and leverage standard compressed sensing theory to derive recovery guarantees.

Unfortunately, this approach fails to deliver optimal rates. The specific reason is that there are $s$-term polynomials that require approximately $m = \ord{s^2}$ i.i.d.\ samples to be stably recovered. The higher level reason is that the assumed \textit{low-complexity} model (functions in $\cH(p)$ have approximately sparse coefficients) is simply too crude. To lower the sample complexity requirement, we refine the model. We do this via \textit{weighted sparsity} \cite{rauhut2016interpolation}, which encodes the additional information that coefficients are both approximately sparse and \textit{decaying}.

\begin{exampleresult}[Near-optimal learning via polynomials, Theorem \ref{thm:main-res-poly}]
\textit{There is a method (based on weighted $\ell^1$-minimization) for learning functions from the data \ef{training-data-intro} that achieves an $L^2_{\varrho}$-norm error that decays like 
\be{
\label{poly-learn-rate}
\ord{(m/\log^4(m))^{1/2-1/p}},\quad m \rightarrow \infty,
}
with high probability, for functions in $\cH(p,\mathsf{M})$ and any $0 < p < 1$.}
\end{exampleresult}

Comparing with the lower bound of $m^{1/2-1/p}$, we see that this procedure is optimal up to the logarithmic term. Note that in the full theorem, we also account for other errors in learning process: namely, \textit{measurement error} (i.e., the effect of the terms $e_i$ in \ef{training-data-intro}), \textit{physical discretization error} (i.e., the effect of discretizing the space $\cV$ in order to perform computations) and \textit{optimization error} (i.e., inexact solution of the optimization problem).

\subsubsection*{DNN existence theory (\S \ref{s:dnn-existence})}

As we discuss further in \S \ref{ss:literature}, the last five years have seen a growing interest in the application of Deep Learning (DL) to challenging parametric model problems. DL has the potential to achieve significant performance gains. Yet, it is poorly understood from a theoretical perspective, especially in terms of the sample complexity, i.e., the amount of training data need to learn good DNN approximations for specific classes of functions.

We commence by briefly reviewing the approximation theory of DNNs, which has been an area of significant research in the last few years. Broadly, this area aims to establish \textit{existence theorems}. Building on the classical \textit{universal approximation theory} of shallow Neural Networks (NNs), these results assert the existence of DNNs of a given complexity -- in terms of the width and depth, number of nonzero weights and biases, or other pertinent metrics -- that approximate functions in a given class to within a prescribed accuracy.

Such results are typically proved through \textit{emulation}. One first shows that DNNs can emulate a standard approximation scheme (e.g., polynomials, piecewise polynomials, splines, wavelets, and so forth), either exactly or up to some error that can be made arbitrarily small. Then one leverages known bounds for the standard scheme to derive the corresponding existence theorem.

Having briefly reviewed this literature, we then describe an exemplar existence theorem for holomorphic functions, which is obtained by emulating a certain $s$-term polynomial approximation.

\begin{exampleresult}[Existence theorem for DNNs, Theorem \ref{thm:existence-holo}]
\textit{Let $(\cV,\nm{\cdot}_{\cV}) = (\bbR,\abs{\cdot})$. Then for any $\bm{b} \in [0,\infty)^{\bbN}$ with $\bm{b} \in \ell^p(\bbN)$ for some $0 < p <1$ and $s \in \bbN$, there is a family $\cN$ of DNNs of width $\ord{s^2}$ and depth $\ord{\log(s)}$ with the following property. For any $f \in \cH(\bm{b})$ there is an element of $\cN$ that achieves an $L^2_{\varrho}$-norm error of $\ord{s^{1/2-1/p}}$.}
\end{exampleresult}

\subsubsection*{Practical existence theory: near-optimal DL (\S \ref{s:learning-proof})}

While important, existence theorems such as this say little about the performance of DL, i.e., computing a DNN from training data by  minimizing a loss function. As noted, analysis of this practical scenario is currently lacking, certainly in terms of sample complexity. It has also been well documented that there is a gap between existence theory and the practical performance of DNNs when trained from finite numbers of samples \cite{abdeljawad2023sampling,adcock2021gap,grohs2023proof}. 

In this section, we describe a theoretical framework that seeks to narrow this gap, termed \textit{practical existence theory} \cite{adcock2021deep,adcock2023near,franco2024practical}. The aim is to show that there not only exists a family DNNs, but also a training strategy akin to standard practice -- i.e., minimizing a loss function -- that attains the near-optimal rates.

\begin{exampleresult}[Practical existence theorem for DNNs, Theorem \ref{thm:main-res-dnn}]
\textit{There is a class of DNNs $\cN$ (of width and depth depending on $m$), a regularized $\ell^2$-loss function and a choice regularization parameter (depending on $m$ only) such that any DNN obtained by solving the resulting training problem (minimizing the loss function over $\cN$) achieves the error decay rate \ef{poly-learn-rate}.}
\end{exampleresult}

Like existence theorems, this result is proved by emulating the sparse polynomial method with DNNs. The class $\cN$ is `handcrafted' in the sense that the weights and biases in the hidden layers are fixed, and chosen specifically to emulate certain polynomials. Only the weights in the final layer are trained.

\subsubsection*{Epilogue: the benefits of practical existence theory and the gap between theory and practice (\S \ref{s:benefits-gap})}

Practical existence theory does not explain the performance of standard DL strategies based on \textit{fully trained models} (i.e., those where all layers of a DNN are trained). Nonetheless, it provides a number of key insights and benefits, which we now list (see \S \ref{s:benefits-gap} for details).

\begin{enumerate}
\item It narrows the gap between theory and practice, by showing the existence of good training procedures, which are near-optimal for certain function classes. 
\item It thereby emphasizes the potential to achieve even better performance in practice with suitably chosen architectures and training strategies.
\item It also stresses the sample complexity aspect of DL for high-dimensional approximation. As argued, this is vital to parametric modelling applications, where data scarcity is of primary concern.
\item The architectures in Theorem \ref{thm:main-res-dnn} are much wider than they are deep. Smoother activation functions also lead to narrower and shallow DNNs. Both insights agree with the empirical performance of fully trained models.
\item Once trained, the DNN in Theorem \ref{thm:main-res-dnn} can be sparsified via \textit{pruning}. This lends credence to the idea that sparse DNNs can perform well in practice.
\item The methodology of practical existence theory is flexible, and can be applied to other problems, such as reduced order models based on deep autoencoders \cite{franco2024practical} and Physics-Informed Neural Networks (PINNs) for PDEs \cite{brugiapaglia2024physics-informed}.
\end{enumerate}

While narrowed, the gap between theory and practice still persists. With an eye towards future research, we end \S \ref{s:benefits-gap} with some further discussion. We discuss how current techniques for analyzing fully trained models may at best achieve a rate of $\ord{m^{-1/2}}$, regardless of $p$, which is strictly slower than those asserted by practical existence theorems. However, the far greater expressivity of DNNs -- in particular, their ability to approximate continuous and discontinuous functions alike -- sets them apart from more standard scientific computing tools such as polynomials. While the latter are near-optimal for holomorphic functions, the flexibility of former is a significant potential advantage.

\subsection{Further literature}\label{ss:literature}

See \cite{2017handbook,le-maitre2010spectral,smith2013uncertainty,sullivan2015introduction} for general introductions to computational UQ and \cite{adcock2022sparse,cohen2015approximation} for more on parametric models and parametric DEs. Note that in the context of UQ, orthogonal polynomial expansions are often termed polynomial \textit{chaos} expansions (see, e.g., \cite{ghanem2003stochastic}).
Best $s$-term polynomial approximation is a type of nonlinear approximation \cite{devore1998nonlinear}. It was developed in a series of works in the context of parametric DEs \cite{beck2014convergence,beck2012optimal,bieri2010sparse,bonito2021polynomial,chkifa2015breaking,cohen2010convergence,cohen2011analytic,hansen2013analytic,todor2007convergence,tran2017analysis}, focusing on Taylor, Legendre or Chebyshev polynomial expansions in infinite dimensions. See, e.g., \cite{cohen2015approximation} and \cite[Chpt.\ 3]{adcock2022sparse} for in-depth reviews. See \S \ref{ss:holo-pdes} for a discussion holomorphic regularity of solutions of various parametric DEs, including the $(\bm{b},\varepsilon)$-holomorphic class and relevant references. Our results on limits of learnability are from \cite{adcock2024optimal} and, as noted, use ideas from information-based complexity. For overviews of this topic, see \cite{novak1988deterministic,novak2008tractability,novak2010tractability,traub1988information}.

In tandem with best $s$-term polynomial approximation theory, there has been a focus on learning polynomial approximations from data. Some early approaches in the context of parametric DEs included (adaptive) interpolation using sparse grids \cite{back2011stochastic,babuska2007stochastic,chkifa2014high-dimensional,ganapathysubramanian2007sparse,gunzburger2014stochastic,ma2009adaptive,mathelin2005stochastic,nobile2008anisotropic,nobile2008sparse,xiu2005high-order}. Another significant line of investigation has been on least-squares methods \cite{chkifa2015discrete,cohen2013stability,migliorati2013polynomial,migliorati2014analysis,berveiller2006stochastic,migliorati2013approximation}. See \cite{guo2020constructing,hadigol2018least} and \cite[Chpt.\ 5]{adcock2022sparse} for reviews.
While simpler than compressed sensing methods, these methods require \textit{a priori} knowledge of a good multi-index set in which to construct the polynomial approximation, and are therefore best suited to the simpler \textit{known anisotropy} setting. \textit{Adaptive} least-squares methods \cite{migliorati2015adaptive,migliorati2019adaptive} can address unknown anisotropy, but they currently lack theoretical guarantees. Compressed sensing for learning polynomial approximations began with the works of \cite{blatman2011adaptive,doostan2011non-adapted,mathelin2012compressed,rauhut2012sparse,yan2012stochastic,yang2013reweighted}, and has seen much subsequent development. See \cite{hampton2017compressive} and \cite[Chpt.\ 7]{adcock2022sparse} and references therein. The focus in this work on \textit{weighted} sparsity originated in \cite{rauhut2016interpolation}. See also \cite{adcock2017infinite-dimensional,adcock2019correcting,peng2014weighted} and references therein. The results described in this work are based on \cite{adcock2024efficient}.

Early works on the application of DL to parametric DEs include \cite{adcock2021deep,cyr2020robust,dal-santo2020data,geist2021numerical,khoo2021solving,laakmann2021efficient}. See also \cite{becker2023learning,cicci2022deep-hyromnet,heiss2021neural,heiss2023multilevel,khara2021neufenet:,lei2022solving,scarabasio2022deep} and references therein  for more recent developments. Closely related to this topic is that of \textit{operator learning} with DNNs -- an area of significant current interest \cite{bhattacharya2021model,boulle2023mathematical,kovachki2023neural,kovachki2024operator,li2021fourier,lu2021learning}. 
One can view holomorphic, Banach-valued functions as an example of holomorphic operators with specific parameterizations of the input space \cite{herrmann2022neural,schwab2023deep}. These operators constitute an important example in the field of operator learning for which there are strong guarantees on both parametric complexity (i.e., the number of DNN parameters required) and sample complexity. The relevant theory for these examples is based on the theory of holomorphic, Banach-valued functions described in this work. See \cite[\S 5.2]{kovachki2024operator} or \cite[\S 3.4]{lanthaler2023operator} for a discussion.

See \cite{cybenko1989approximation,hornik1989multilayer} and \cite{pinkus1999approximation} for the classical \textit{universal approximation theory} of (shallow) NNs. The modern study of existence theory was initiated in \cite{yarotsky2017error}, and has since seen a wealth of developments. We review this literature in \S \ref{s:dnn-existence}. The polynomial emulation results used in this work are based on \cite{daws2019analysis,de-ryck2021approximation,opschoor2022exponential}. 

The theory-to-practice gap in DL was studied empirically \cite{adcock2021gap} and theoretically in \cite{abdeljawad2023sampling,grohs2023proof}. Practical existence theorems were first established in \cite{adcock2021deep,adcock2021gap}. Our results are based on \cite{adcock2023near}. This approach is not unique to high-dimensional regression problems. Similar ideas have been used to assert that it is possible to compute stable, accurate and efficient DNNs for inverse problems in imaging \cite{colbrook2022difficulty,neyra-nesterenko2023stable}. Here, one starts with a standard regularization problem, such as TV minimization, then exploits the fact that $n$ steps of a first-order optimization method for solving this problem can be reinterpreted as a DNN with $\ord{n}$ layers -- a process known as \textit{unrolling} \cite{monga2021algorithm} (see also \cite[Chpts.\ 19-21]{adcock2021compressive}). Interestingly, this also leads to principled ways to design DNN architectures for DL in inverse problems \cite{monga2021algorithm}, although the trained networks may not be robust \cite{antun2023ai,antun2020instabilities}.

\section{Problem statement and notation}\label{s:prob-statement}

In this section, we formalize the problem studied in this work and establish important notation.

Throughout this work, $\cU = [-1,1]^{\bbN}$ and $(\cV,\nm{\cdot}_{\cV})$ is a Banach space over $\bbR$. We write $\bm{y} = (y_i)_{i \in \bbN}$ for the independent variable in $\cU$. 
We equip $\cU$ with the uniform probability measure $\varrho$ and, for $1 \leq p \leq \infty$, write $L^p_{\varrho}(\cU;\cV)$ for the Lebesgue--Bochner space of (equivalence classes of) strongly $\varrho$-measurable functions $f: \cU \rightarrow \cV$ for which $\nm{f}_{L^{p}_{\varrho}(\cU ; \cV)} < \infty$, where 
\be{
\nm{f}_{L^p_{\varrho}(\cU;\cV)} : = 
\begin{cases} 
\left( \int_{\cU} \nm{f( \bm{y})}_{\cV}^p \D \varrho (\bm{y}) \right)^{1/p} & 1 \leq p < \infty ,
\\
\mathrm{ess} \sup_{\bm{y} \in \cU} \nm{f(\bm{y})}_{\cV}  & p = \infty .
\end{cases}
\label{L_p_U_V}
}
When $\cV = (\bbR,\abs{\cdot})$, we just write $L^{p}_{\varrho}(\cU)$ for the corresponding Lebesgue space of real-valued functions.

Let $f \in L^2_{\varrho}(\cU ; \cV)$ be the unknown target function and consider sample points $\bm{y}_1,\ldots,\bm{y}_m \sim_{\mathrm{i.i.d.}} \varrho$. Then we consider training data 
\be{
\label{training-data}
\{ (\bm{y}_i,f(\bm{y}_i) + e_i ) \}^{m}_{i=1} \subseteq \cU \times \cV,
}
where $e_i \in \cV$ represents measurement error. In this work, we consider adversarial noise of small norm, i.e., $\sum^{m}_{i=1} \nm{e_i}^2_{\cV} \ll 1$. Statistical models can also be considered -- see, e.g., \cite{migliorati2015convergence}. The problem is then to learn an approximation to $f$ based on the data \ef{training-data}.

We require some notation for sequences indexed via possibly multi-indices. Let $d = \bbN \cup \{ \infty \}$. We write $\bm{\nu} = (\nu_k)^{d}_{k=1}$ for an arbitrary multi-index in $\bbN^d_0$. Let $\Lambda \subseteq \bbN^d_0$ be a finite or countable set of multi-indices. For $0 < p \leq \infty$, the space $\ell^p(\Lambda ; \cV)$ consists of all $\cV$-valued sequences $\bm{c} = (c_{\bm{\nu}})_{\bm{\nu} \in \Lambda}$ for which $\nm{\bm{c}}_{p;\cV} < \infty$, where
\eas{
\nm{\bm{c}}_{p;\cV} = \begin{cases} \left ( \sum_{\bm{\nu} \in \Lambda} \nm{c_{\bm{\nu}}}^p_{\cV} \right )^{\frac1p} & 0 < p < \infty ,
\\
\sup_{\bm{\nu} \in \Lambda} \nm{c_{\bm{\nu}}}_{\cV} & p = \infty. 
\end{cases}
}
When $(\cV,\nm{\cdot}_{\cV}) = (\bbR,\abs{\cdot})$, we just write $\ell^p(\Lambda)$ and $\nm{\cdot}_{p}$. Given a sequence $\bm{c} = (c_{\bm{\nu}})_{\bm{\nu} \in \Lambda}$, we define its support
\bes{
\mathrm{supp}(\bm{c}) = \{ \bm{\nu} \in \Lambda : c_{\bm{\nu}} \neq 0 \} \subseteq \Lambda.
}
Now let $\Lambda = \bbN$ and $0 < p \leq \infty$. Given a real-valued sequence $\bm{z} = (z_i)_{i \in \bbN} \in \ell^{\infty}(\bbN)$, we define its minimal monotone majorant
\begin{equation}
\label{min-mon-maj}
\tilde{\bm{z}} = (\tilde{z}_i)_{i \in \bbN},\quad \text{where }
\tilde{z}_i = \sup_{j \geq i} | z_{j}|,\ \forall i \in \bbN
\end{equation}
and the \textit{monotone $\ell^p$-space} $\ell^p_{\mathsf{M}}(\bbN)$ as
\begin{equation}
\label{monotone-lp}
\ell^p_{\mathsf{M}}(\bbN) = \{ \bm{z} \in \ell^{\infty}(\bbN) : \nm{\bm{z}}_{p,\mathsf{M}} : = \nm{\tilde{\bm{z}}}_{p} < \infty  \}.
\end{equation}
Let $d \in \bbN \cup \{ \infty \}$. We write $\bm{e}_j$, $j = 1,\ldots,d$, for the canonical basis vectors in $\bbR^d$.
Given multi-indices $\bm{\nu} = (\nu_k)^{d}_{k=1}$ and $\bm{\mu} = (\mu_k)^{d}_{k=1}$, the inequality $\bm{\nu} \geq \bm{\mu}$ is interpreted componentwise, i.e., $\nu_k \geq \mu_k$, $\forall k$. We define $\bm{\nu} > \bm{\mu}$ analogously. We write $\bm{0}$ for the multi-index of zeros and $\bm{1}$ for the multi-index of ones.

\rem{[Other measures and domains]
We work with the uniform probability measure on $\cU = [-1,1]^{\bbN}$ for convenience. This means that we construct polynomial approximations using the corresponding multivariate Legendre polynomials. It is possible to work more generally, for instance by considering ultraspherical measures and the resulting Jacobi polynomials. See  \cite{adcock2024optimal}. The cases where $\cU = \bbR^{\bbN}$ with the Gaussian measure (corresponding to Hermite polynomials) or $\cU = [0,\infty)^{\bbN}$ with the gamma measure (corresponding to Laguerre) have also been studied, although the theory is less complete. Results on best $s$-term polynomial approximation are known in this case, as are DNN existence theorems \cite{dung2023deep,schwab2023deepa}. However, results akin to those we present in this work on learning sparse polynomial or DNN approximations from data are lacking.
}

\rem{
[Error metric]
Throughout this work, we will measure the error of the various approximations in the $L^2_{\varrho}(\cU ; \cV)$-norm. The majority of the results we present extend to the stronger $L^{\infty}_{\varrho}(\cU ; \cV)$-norm. In particular, the various upper bounds that establish algebraic rates with index $1/2-1/p$ can be modified to show algebraic rates with index $1-1/p$ (see, e.g., \cite{adcock2024efficient}). For succinctness, we do not present this modification.
}

\section{Holomorphic functions of infinitely many variables}\label{s:holo-fns}

In this section, we first formally define the classes of holomorphic functions considered in this work and then discuss their relevance to parametric DEs.

\subsection{$(\bm{b},\varepsilon)$-holomorphic functions}

We commence with the definition of holomorphy.

\begin{definition}
[Holomorphy]
Let $\mathcal{O} \subseteq \mathbb{C}^{\bbN}$ be an open set and  $\cV$ be a Banach space. A  function $f: \mathcal{O} \rightarrow \cV$ is \textit{holomorphic in $\mathcal{O}$} if it is holomorphic with respect to each variable  in  $\mathcal{O}$. That is to say, for any $\bm{z} \in \mathcal{O}$ and any $j \in \bbN$, the following limit exists in $\cV$:
\[
\lim_{\substack{h \in \bbC \\ h \rightarrow 0}}  \dfrac{f(\bm{z}+h \bm{e}_j)-f(\bm{z})}{h} \in \cV.
\]
\end{definition}

\defn{
[Holomorphic extension]
Let $\cU \subset \cO \subseteq \bbC^{\bbN}$ be an open set. A function $f : \cU \rightarrow \cV$ has a \textit{holomorphic extension to $\cO$} (or simply, \textit{is holomorphic in $\cO$}) if there is a $\tilde{f} : \cO \rightarrow \cV$ that is holomorphic in $\cO$ for which $\tilde{f} |_{\cU} = f$.  In this case, we also define $\nm{f}_{L^{\infty}(\cO; \cV)}:=\nm{{\tilde{f}}}_{L^{\infty}(\cO; \cV)} = \sup_{\bm{z} \in \cO} \nm{\tilde{f}(\bm{z})}_{\cV}$ or, when $\cV = \bbC$, simply $\nm{f}_{L^{\infty}(\cO)}$. If $\cO$ is a closed set, then we say that $f$ is holomorphic in $\cO$ if it has a holomorphic extension to some open neighbourhood of $\cO$.
}

In this work, we consider functions that possess holomorphic extensions to (unions of) tensor-products of Bernstein ellipses. The \textit{Bernstein ellipse} with parameter $\rho > 1$ is the set
\bes{
\cE(\rho) = \left\{ (z+z^{-1})/2 : z \in \bbC,\ 1 \leq | z | \leq \rho \right \} \subset \bbC.
}
See, e.g., \cite[Chpt.\ 8]{trefethen2013approximation}.
This is an axis-aligned ellipse with foci at $\pm 1$ and with major and minor semi-axis lengths given by $(\rho + 1/\rho)/2$ and $(\rho-1/\rho)/2$, respectively. By convention, we define $\cE(\rho) = [-1,1]$ when $\rho = 1$. Next, we define the \textit{Bernstein polyellipse} with parameter $\bm{\rho} = (\rho_i)_{i \in \bbN} > \bm{1}$ by
\begin{equation*}
\cE(\bm{\rho}) = \cE(\rho_1) \times \cE(\rho_2) \times \cdots \subset \bbC^{\bbN}.
\end{equation*}
We are now ready to define the class of $(\bm{b},\varepsilon)$-holomorphic functions \cite{chkifa2015breaking,schwab2019deep}.

\defn{
[$(\bm{b},\varepsilon)$-holomorphic functions]
\label{d:b-eps-holomorphy}
Let $\bm{b} = (b_{j})_{j \in \bbN} \in [0,\infty)^{\bbN}$ with $\bm{b} \in \ell^1(\bbN)$ and $\varepsilon > 0$. A function $f : \cU \rightarrow \cV$ is \textit{$(\bm{b},\varepsilon)$-holomorphic} if it has a holomorphic extension to the region
\begin{equation}\label{def_R_b_e}
\cR({\bm{b},\varepsilon}) =  \bigcup \left\lbrace \cE(\bm{\rho}): \bm{\rho} \geq \bm{1},\  \sum_{j=1}^{\infty} \left( \dfrac{\rho_j+\rho_j^{-1}}{2} -1 \right) b_j \leq  \varepsilon \right\rbrace \subset \bbC^{\bbN}.
\end{equation}
}
We discuss the motivations for this definition next. Beforehand, it is worth noting that the parameter $\varepsilon$ is technically redundant, since we can always rescale $\bm{b}$. In the remainder of this work, we assume that $\varepsilon = 1$. For later use, we now define the class of $(\bm{b},1)$-holomorphic functions with norm at most one as
\bes{
\cH(\bm{b}) = \left \{ f : \cU \rightarrow \cV\text{ $(\bm{b},1)$-holomorphic} : \nm{f}_{L^{\infty}(\cR(\bm{b},1);\cV)}  \leq 1 \right \}.
}

\rem{[Functions of finitely many variables]
\label{rem:fin-dim}
In finite dimensions (in particular, $d = 1$), one normally considers functions that are holomorphic in a single Bernstein polyellipse $\cE(\bm{\rho})$. For such functions, it is well known that one can find polynomial approximations that converge exponentially fast. As we see later, in infinite dimensions we use the assumption of holomorphy in the region \ef{def_R_b_e} in order to obtain algebraic rates of convergence.

It is, however, worth noting that finite-dimensional holomorphy in a Bernstein polyellipse implies $(\bm{b},\varepsilon)$-holomorphy. Let $f : [-1,1]^d \rightarrow \bbR$ be a function of finitely many variables that is holomorphic in $\cE(\bar\rho_1) \times \cdots \times \cE(\bar\rho_d) \subset \bbC^d$. We can extend $f$ to a function $\bar{f} : \cU \rightarrow \bbC$ in the standard way as $\bar{f}(y_1,y_2,\ldots) = f(y_1,\ldots,y_d)$, $\forall \bm{y} = (y_i)_{i \in \bbN} \in \cU$. Now define
\eas{
b_i = \varepsilon ( (\bar\rho_i + \bar\rho^{-1}_{i})/2 - 1  )^{-1},\ i \in [d],\qquad b_i = 0,\ i \in \bbN \backslash [d].
}
Then $\bar{f}$ is $(\bm{b},\varepsilon)$-holomorphic. Therefore, all the results that follow on approximation of  $(\bm{b},\varepsilon)$-holomorphic functions also apply \textit{mutatis mutandis} to functions of $d$ variables that are holomorphic in a single Bernstein ellipse.
}

\subsection{Holomorphy and parametric DEs}\label{ss:holo-pdes}

As mentioned, $(\bm{b},\varepsilon)$-holomorphic functions were first developed in the context of parametric DEs. There is now a wealth of literature that establishes that solution maps of many parametric DEs posses this regularity.

The classical example of a parametric PDEs is the parametric stationary diffusion equation with parametrized diffusion coefficient
\be{
\label{eq:fruitfly}
-\nabla_{\bm{x}} \cdot (a(\bm{x},\bm{y}) \nabla_{\bm{x}} u(\bm{x}, \bm{y})) = F(\bm{x}), \  \bm{x} \in \Omega,\qquad u(\bm{x},\bm{y}) = 0, \ \bm{x} \in \partial \Omega.
}
Here $\bm{x} \in \Omega $ is the spatial variable, $\Omega \subset \mathbb{R}^k$ is the domain domain, which is assumed to have Lipschitz boundary, and $\nabla_{\bm{x}}$ is the gradient with respect to $\bm{x}$. The function $F$ is the forcing term, and is assumed to be nonparametric.
We assume that \eqref{eq:fruitfly} satisfies a \textit{uniform ellipticity condition}
\be{
\label{uniform-ellipticity}
 \essinf_{\bm{x} \in \Omega} a(\bm{x},\bm{y}) \geq r,\quad \forall \bm{y} \in \cU,
}
for some $r > 0$. Now let $\cV = H^1_0(\Omega)$ be the standard Sobolev space of functions with weak first-order derivatives in $L^2(\Omega)$ and traces vanishing on $\partial \Omega$. Then the problem \eqref{eq:fruitfly} has a well-defined \textit{parametric solution map}
\be{
\label{u-soln-map}
u : \cU \rightarrow \cV,\ \bm{y} \mapsto u(\cdot,\bm{y}).
}
Now suppose that the diffusion coefficient has the affine parametrization
\be{
\label{eq:affine-diffusion}
a(\bm{x},\bm{y}) = a_0(\bm{x}) + \sum^{\infty}_{j=1} y_j \psi_j(\bm{x}),
}
where the functions $a_0, \psi_1,\psi_2,\ldots \in L^{\infty}(\Omega)$. Notice that uniform ellipticity \eqref{uniform-ellipticity} for this problem is equivalent to the condition
\bes{
\sum^{\infty}_{j=1} | \psi_j(\bm{x}) | \leq a_0(\bm{x}) - r,\quad \forall \bm{x} \in \Omega.
}
Under this assumption, one can show (see, e.g., \cite[Prop.\ 4.9]{adcock2022sparse}.) that the solution map \eqref{u-soln-map} is well defined and $(\bm{b},\varepsilon)$-holomorphic for $0 < \varepsilon < r$ with
\be{
\label{b-affine}
\bm{b} = (b_i)^{\infty}_{i=1},\quad \text{with }b_i = \nm{\psi_i}_{L^{\infty}(\Omega)},\ \forall i \in \bbN.
}
This example was first considered in \cite{bieri2010sparse,cohen2010convergence,cohen2011analytic}. Subsequently, these results were generalized to many other classes of parametric DEs and PDEs. One such extension (see, e.g., \cite{chkifa2015breaking}) considers the problem \ef{eq:fruitfly} with various nonaffine parametric diffusion coefficients, such as 
\bes{
a(\bm{x},\bm{y}) =a_0(\bm{x}) + \left ( \sum^{\infty}_{j=1} y_j \psi_j(\bm{x}) \right )^2\quad \text{or}\quad a(\bm{x},\bm{y}) =\exp\left ( \sum^{\infty}_{j=1} y_j \psi_j(\bm{x}) \right ).
}
This is part of a general framework for showing holomorphy for parametric weak problems in Hilbert spaces. See \cite[Thm.\ 4.1]{chkifa2015breaking}, \cite[Cor.\ 2.4]{cohen2015approximation} or \cite[\S 4.3.1]{adcock2022sparse}. See also \cite{kunoth2013analytic,rauhut2017compressive}. Another related general framework considers parametric implicit operator equations (see \cite[Thm.\ 4.3]{chkifa2015breaking}, \cite[Thm.\ 2.5]{cohen2015approximation} or \cite[\S 4.3.1]{adcock2022sparse}).

These frameworks can be used to establish holomorphy results for various problems  beyond \ef{eq:fruitfly}. This includes parabolic problems (see \cite[\S 5.1]{chkifa2015breaking} or \cite[\S 2.2]{cohen2015approximation}) and various types of nonlinear, elliptic PDEs (see \cite[\S 5.2]{chkifa2015breaking} or \cite[\S 2.3]{cohen2015approximation}). Another prominent example includes PDEs such as \ef{eq:fruitfly} over parametrized domains $\Omega = \Omega_{\bm{y}}$. Holomorphy of the solution map is often referred to as \textit{shape holomorphy} in this context. See \cite[\S 5.3]{chkifa2015breaking}, \cite[\S 2.2]{cohen2015approximation}, \cite{castrillon2016analytic,cohen2018shape} and references therein. Other examples include parametric hyperbolic problems \cite{hoang2012regularity} and certain classes of parametric control problems \cite{kunoth2013analytic}.

For further reviews, see \cite{cohen2015approximation} or \cite[Chpt.\ 4]{adcock2022sparse}. Note that the above problems all involve PDEs. Classes of parametric ODEs also admit holomorphic regularity. See \cite[\S 4.1]{adcock2022sparse} and \cite{hansen2013sparse}. Recent work in this direction has also shown such regularity for parametric ODEs arising from diffusion on graphs \cite{ajavon2024surrogate}.

\subsection{Known and unknown anisotropy}

Functions that are $(\bm{b},\varepsilon)$-holomorphic are \textit{anisotropic}: they depend more smoothly on some variables rather than others. This behaviour is dictated by the parameter $\bm{b}$, with a larger value of an entry $b_i$ corresponding to less smoothness with respect to the variable $y_i$. This can be seen by noting that the condition
\bes{
\sum^{\infty}_{i=1} \left ( \frac{\rho_i + \rho^{-1}_i}{2} - 1 \right ) b_i \leq 1,
}
in \ef{def_R_b_e} holds only for smaller values of $\rho_i$ when $b_i$ is large, meaning that $f$ only admits an extension in the $y_i$ variable to a relatively small Bernstein ellipse.

In practice, \textit{a priori} analysis of a given problem may establish that the target function is $(\bm{b},\varepsilon)$-holomorphic for some $\bm{b}$. However, an optimal value for $\bm{b}$ may be unknown. This situation arises in parametric DEs. For problems such as the stationary diffusion equation \ef{eq:fruitfly} with affine diffusion, one can find a sufficient value \ef{b-affine} of $\bm{b}$ for which the parametric solution map is $(\bm{b},\varepsilon)$-holomorphic, although this value may not be sharp.  In general, since $\cR(\bm{b},\varepsilon) \subseteq \cR(\bm{b}',\varepsilon)$ whenever $\bm{b}' \leq \bm{b}$, it is difficult to know whether a value of $\bm{b}$ obtained from some analysis is optimal. Further, for more complicated problems, such as some of those described above, one can derive theoretical guarantees that assert holomorphy, but without an estimate for the region itself \cite[Rem.\ 2.6]{cohen2015approximation}.

For this reason, the primary focus of this work is the \textit{unknown anisotropy} setting. We may assume the target function $f$ is $(\bm{b},\varepsilon)$-holomorphic for some $\bm{b}$, but we do not have access to $\bm{b}$ itself. In particular, this means that we cannot use $\bm{b}$ to design a method for learning $f$ from data -- the holomorphy assumption can only be used to provide bounds for the resulting approximation error. 

\rem{
To emphasize this consideration further, consider the function
\bes{
f(y_1,y_2) = \sin(1000 y_2) / (1.1-y_1).
}
This function is entire with respect to the variable $y_2$, but only has a pole at $y_1 = 1.1$. By Remark \ref{rem:fin-dim}, it is $(\bm{b},1)$-holomorphic for $\bm{b} = (b_1,b_2,\ldots) = (10,0,0,\ldots)$. When building a polynomial approximation to $f$, knowledge of $\bm{b}$ may lead one to include only low-degree terms in $y_2$ and more higher-degree terms in $y_1$. Yet this is completely the opposite of what one should do in practice. The function $\sin(1000 y_2)$, while entire, is highly oscillatory, and can only be resolved by using high-degree polynomials. On the other hand, the function
\bes{
g(y_1,y_2) = 1/(1.1-y_1)
}
is also $(\bm{b},1)$-holomorphic with the same $\bm{b}$. Yet, it requires no non-constant terms in $y_2$ in order to approximate it accurately with a polynomial.
}

\subsection{$\ell^p$-summability and the $\cH(p)$ and $\cH(p,\mathsf{M})$ classes}\label{ss:holo-p-classes}

As we will see in \S \ref{s:sparse-poly}, it is generally impossible to approximate a $(\bm{b},\varepsilon)$-holomorphic function without a further assumption on the sequence $\bm{b}$. In particular, the terms of $\bm{b}$ need to decay sufficiently fast. Henceforth, we will assume that $\bm{b}$ is not just in $\ell^1(\bbN)$, but that,  for some $0 < p <1$, either 
\bes{
\bm{b} \in \ell^p(\bbN)\text{ or }\bm{b} \in \ell^p_{\mathsf{M}}(\bbN).
}
Here we recall that $\ell^p_{\mathsf{M}}(\bbN)$ is the monotone $\ell^p$-space \ef{monotone-lp}. When considering the unknown anisotropy setting, we also define the function classes
\be{
\cH(p) =  \bigcup \left \{ \cH(\bm{b}) : \bm{b} \in \ell^p(\bbN),\ \bm{b} \in [0,\infty)^{\bbN},\   \nm{\bm{b}}_p \leq 1 \right \} 
}
and
\be{
\label{HpM}
\cH(p,{\mathsf{M}})  = \bigcup  \left \{ \cH(\bm{b}) : \bm{b} \in \ell^p_{\mathsf{M}}(\bbN),\ \bm{b} \in [0,\infty)^{\bbN},\  \nm{\bm{b}}_{p,\mathsf{M}} \leq 1 \right \} .
}
In what follows we will ask for a method to provide a uniform error bound (depending on $m$, the amount of data, and $p$) for \textit{any} function in these classes.

\section{Best $s$-term polynomial approximation}\label{s:sparse-poly}

In this section, we make a first foray into polynomial approximation of holomorphic functions. We first introduce 
orthogonal polynomials and orthogonal polynomial expansions in $L^2_{\varrho}(\cU;\cV)$. We then introduce best $s$-term polynomial approximation and provide a theorem on algebraic convergence for it. See, e.g., \cite[Chpts. 2-3]{adcock2022sparse} or \cite{chkifa2015breaking,cohen2015approximation}, for further information on this material.

\subsection{Orthogonal polynomials}

Let $P_0,P_1,\ldots$ denote the classical one-dimensional Legendre polynomials
 with the normalization
\be{
\label{P-max}
\nm{P_{\nu}}_{L^{\infty}([-1,1])} = P_{\nu}(1) = 1.
}
We consider the orthonormalized version of these polynomials. Since
\bes{
\int^{1}_{-1} | P_{\nu}(y) |^2 \D y = \frac{2}{2 \nu+1},\quad \forall \nu \in \bbN_0,
}
we define these as 
\be{
\label{psi-def}
\psi_{\nu}(y) = \sqrt{2 \nu+1} P_{\nu}(y),\quad \forall \nu \in \bbN_0.
}
The set $\{ \psi_{\nu} \}_{\nu \in \bbN_0} \subset L^2_{\rho}([-1,1])$ forms an orthonormal basis, where $\rho$ is the one-dimensional uniform measure.

We construct an orthonormal basis of $L^2_{\varrho}(\cU)$ by tensorization. Let 
\bes{
\cF = \left \{ \bm{\nu} = (\nu_k )^{\infty}_{k=1} \in \bbN^{\bbN}_0 : | \{ k : \nu_k \neq 0 \} | < \infty \right \}
}
be the set of multi-indices with finitely many nonzero terms and set
\bes{
\Psi_{\bm{\nu}}(\bm{y}) = \prod_{i \in \bbN} \psi_{\nu_i}(y_i),\quad \forall \bm{y} \in \cU,\ \bm{\nu} \in \cF.
}
Note that $\psi_{0} = 1$ by construction. Therefore, this is equivalent to a product over finitely many terms, i.e.,
\be{
\label{Psi-def}
\Psi_{\bm{\nu}}(\bm{y}) = \prod_{i : \nu_i \neq 0} \psi_{\nu_i}(y_i).
}
Given these functions, it can now be shown that the set
\be{
\label{ONB}
\{ \Psi_{\bm{\nu}} \}_{\bm{\nu} \in \cF} \subset L^2_{\varrho}(\cU)
}
constitutes an orthonormal basis for $L^2_{\varrho}(\cU)$ \cite[\S 3]{cohen2015approximation}.

Using \ef{P-max}--\ef{Psi-def}, we also deduce that
\be{
\label{u-def}
\nm{\Psi_{\bm{\nu}}}_{L^{\infty}_{\varrho}(\cU ; \cV)} = | \Psi_{\bm{\nu}}(\bm{1}) | = \prod_{k \in \bbN} \sqrt{2 \nu_k + 1} = : u_{\bm{\nu}}.
}
The values $u_{\bm{\nu}}$ will be of use later. For convenience, we also define
\be{
\label{intrinsic-weights}
\bm{u} = (u_{\bm{\nu}})_{\bm{\nu} \in \cF}.
}

\subsection{Orthogonal polynomial expansions}

Let $f \in L^2_{\varrho}(\cU ; \cV)$. Then we have the convergent expansion (in $L^2_{\varrho}(\cU ; \cV)$)
\be{
\label{f_exp}
f = \sum_{\bm{\nu} \in \cF} c_{\bm{\nu}} \Psi_{\bm{\nu}},\quad \text{where }c_{\bm{\nu}} =\int_{\cU} f(\bm{y}) \Psi_{\bm{\nu}}(\bm{y})    \D \varrho (\bm{y}) \in \cV.
}
Note that the \textit{coefficients} $c_{\bm{\nu}}$ are elements of $\cV$ and defined by Bochner integrals. For convenience, we denote the infinite sequence of coefficients of $f$ by
\be{
\label{f-coeff}
\bm{c} = (c_{\bm{\nu}})_{\bm{\nu} \in \cF}.
}
Parseval's identity gives that $\bm{c} \in \ell^2(\cF ; \cV)$ with
$\nm{\bm{c}}_{2;\cV} = \nm{f}_{L^2_{\varrho}(\cU ; \cV)}$.

\subsection{Best $s$-term polynomial approximation}

Let $s \in \bbN$. An \textit{$s$-term polynomial approximation} to $f$ is an approximation
\be{
\label{f-S-def}
f \approx f_{S} = \sum_{\bm{\nu} \in S} c_{\bm{\nu}} \Psi_{\bm{\nu}}
}
for some multi-index set $S \subset \cF$, $|S| = s$. This raises the question: which index set $S$ should one choose? One answer is to select a set $S$ which provides the best approximation, an approach known as \textit{best $s$-term approximation} and itself a type of nonlinear approximation \cite{devore1998nonlinear}. Formally, a \textit{best $s$-term approximation} $f_s$ of $f$ (with respect to the $L^2_{\varrho}$-norm) is defined as
\be{
\label{best-s-term-inf}
f_s = f_{S^*},\qquad \text{where }S^* \in \argmin{} \{ \nm{f - f_S}_{L^2_{\varrho}(\cU ; \cV)} : S \subset \cF,\ |S| = s \}.
}
This approximation can also be characterized in terms of the coefficients of $f$. Due to Parseval's identity, the error of the approximation $f_S$ is precisely
\be{
\label{Parseval}
\nm{f - f_S}^2_{L^2_{\varrho}(\cU;\cV)} = \sum_{\bm{\nu} \in \cF \backslash S} \nm{ c_{\bm{\nu}} }^2_{\cV}.
}
Therefore, any set $S^*$ that yields a best $s$-term approximation consists of multi-indices corresponding to the largest $s$ coefficient norms $\nm{ c_{\bm{\nu}} }_{\cV}$. Specifically,
\bes{
S^* = \{ \bm{\nu}_1,\bm{\nu}_2,\ldots,\bm{\nu}_s \}, 
}
where $\bm{\nu}_1,\bm{\nu}_2,\ldots$ is an ordering of the multi-index set $\cF$ such that $\nm{c_{\bm{\nu}_1} }_{\cV} \geq \nm{ c_{\bm{\nu}_2} }_{\cV} \geq \cdots $. It follows immediately from this and \ef{Parseval} that the error
\be{
\label{best-s-term-error}
\nm{f - f_s}^2_{L^2_{\varrho}(\cU;\cV)} = \sum_{i > s} \nm{c_{\bm{\nu}_i}}^2_{\cV}.
}
is precisely the $\ell^2(\cF;\cV)$-norm of the sequence of coefficients \ef{f-coeff}, excluding those coefficients with indices in $S^*$. Note that $S^*$ (and therefore $f_s$) is generally nonunique. However, this fact causes no difficulties in what follows.

\subsection{Rates of best $s$-term polynomial approximation}

The best $s$-term approximation describes the best possible approximation obtainable with an $s$-term polynomial expansion. It is therefore important to provide bounds for this \textit{benchmark} approximation in 
the case of holomorphic functions. Besides providing insight into limits of approximation, these bounds are also useful when we come to learn sparse polynomial approximations from data.

As we discussed above, Parseval's identity relates questions of best $s$-term approximation to $f$ (in the $L^2_{\varrho}(\cU ; \cV)$ norm) to approximation of its coefficients (in the $\ell^2(\cF ; \cV)$-norm). It motivates one to study best $s$-term approximation of sequences, rather than functions. To this end, we now introduce some additional notation. Let $\Lambda \subseteq \cF$, $0 < p \leq \infty$, $\bm{c} \in \ell^p(\Lambda;\cV)$ and $s\in \mathbb{N}_0$ with $s \leq |\Lambda|$. The \textit{$\ell^p$-norm best $s$-term approximation error} of the sequence $\bm{c}$ is defined as
\be{
\label{best_s_term}
\sigma_{s}(\bm{c})_{p;\cV} = \min \left \{ \nm{\bm{c} - \bm{z}}_{p;\cV} : \bm{z} \in \ell^p(\Lambda;\cV),\ | \supp(\bm{z}) | \leq s \right \}.
}
Note that when $p = 2$ we have
$
 \sigma_{s}(\bm{c})_{2;\cV} = \nm{f - f_s}_{L^2_{\varrho}(\cU ; \cV)},
$
where $f_s$ is a best $s$-term approximation to $f$ \ef{best-s-term-inf}.

We now state a well known result regarding the best $s$-term approximation (of coefficients) of holomorphic functions (see, e.g., \cite[Thm.\ 3.28]{adcock2022sparse} or \cite[\S 3.2]{cohen2015approximation}).
\thm{
[Algebraic convergence of the best $s$-term approximation]
\label{t:best-s-term}
Let $\bm{b} \in [0,\infty)^{\bbN}$ be such that $\bm{b} \in \ell^p(\bbN)$ for some $0 < p < 1$. Then for any $s \in \bbN$ and $p \leq q \leq 2$, there exists a set $S \subset \cF$ with $|S| \leq s$ such that
\bes{
\sigma_{s}(\bm{c})_{p;\cV}  \leq \nm{\bm{c}-\bm{c}_S}_{q;\cV} \leq C(\bm{b},p) \cdot s^{\frac1q-\frac1p},
}
for all $f \in \cH(\bm{b})$ with coefficients $\bm{c}$ as in \ef{f-coeff}.
}

\rem{
[Sharpness of the algebraic rate] 
Theorem \ref{t:best-s-term} provides an upper bound for the best $s$-term approximation error. However, it is also possible to provide a lower bound. As shown in \cite[Thm.\ 5.6]{adcock2023monte}, the are choices of $\bm{b} \in \ell^p(\bbN)$ and functions $f \in \cH(\bm{b})$ for which
\bes{
\limsup_{s \rightarrow \infty} \left \{ s^{\frac1r-\frac12} \cdot \nm{f - f_s}_{L^2_{\varrho}(\cU ; \cV)} \right \}  = + \infty
}
for any $0 < r < p$. Thus, the algebraic exponent $1/2-1/p$ is sharp. 
}

\subsection{How high is high dimensional?}

As in Remark \ref{rem:fin-dim}, let $f : [-1,1]^d \rightarrow \bbR$ be holomorphic in a single Bernstein ellipse $\cE(\bar{\bm{\rho}}) = \cE(\bar\rho_1) \times \cdots \times \cE(\bar\rho_d) \subset \bbC^d$. Then its best $s$-term polynomial approximation converges with exponential rate in $s^{1/d}$: specifically
\bes{
\nm{f - f_s}_{L^{2}_{\varrho}(\cU ; \cV)} \leq \nm{f}_{L^{\infty}(\cR ; \cV)} C(\epsilon) \sqrt{s} \exp \left ( - \left ( \frac{s d! \prod^{d}_{j=1} \log(\bar\rho_j)}{1+\epsilon} \right )^{1/d} \right ),
}
for any $0 < \epsilon < 1$ and all $s \geq \bar{s}$, where $\bar{s}$ is a constant depending on  $d$,$\epsilon$ and $\bar{\bm{\rho}}$ only. See \cite{tran2017analysis} or \cite[Thm.\ 3.21]{adcock2022sparse}.
In low dimensions, this rate accurately predicts the convergence of the best $s$-term approximation. However, for larger $d$ this rate is generally not witnessed for any finite $s$, due to the exceedingly slow growth of the term $s^{1/d}$. For example, $s^{1/d} \leq 3.2$ in $d = 8$ dimensions for all $1 \leq s \leq 10,000$ -- a range of $s$ that is reasonable in practice. 

\begin{figure}[t]
\begin{center}
\begin{small}
 \begin{tabular}{cc}
\includegraphics[width = 0.40\textwidth]{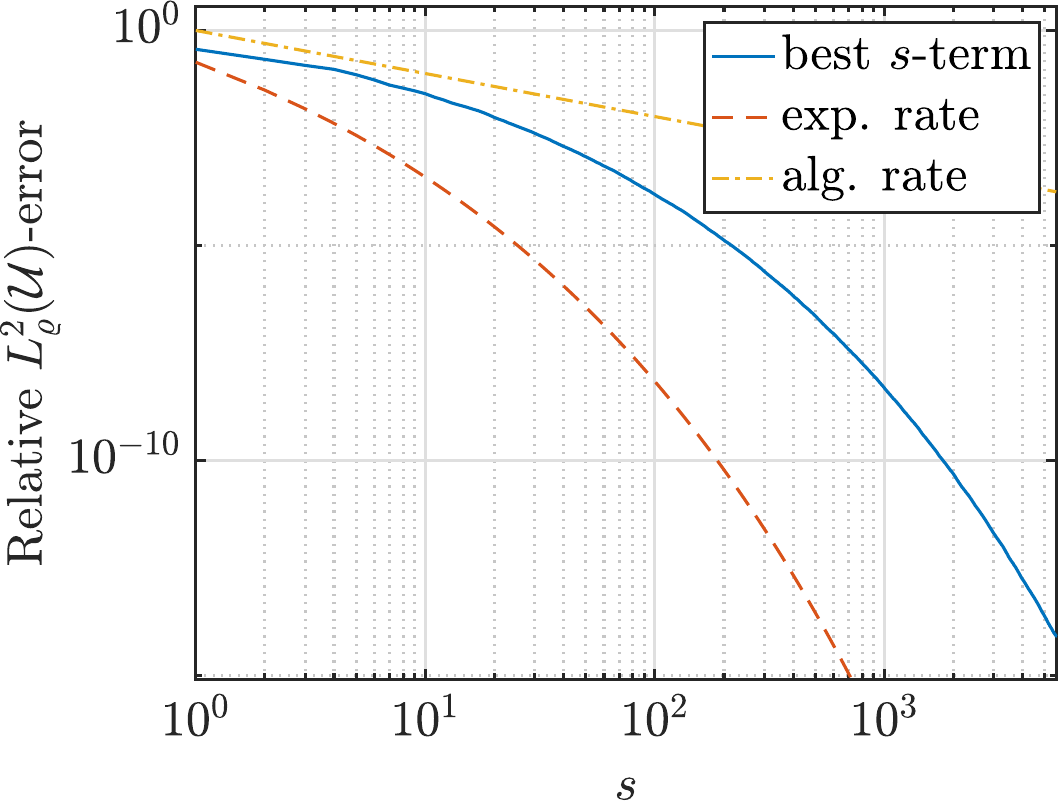}
&
\includegraphics[width = 0.40\textwidth]{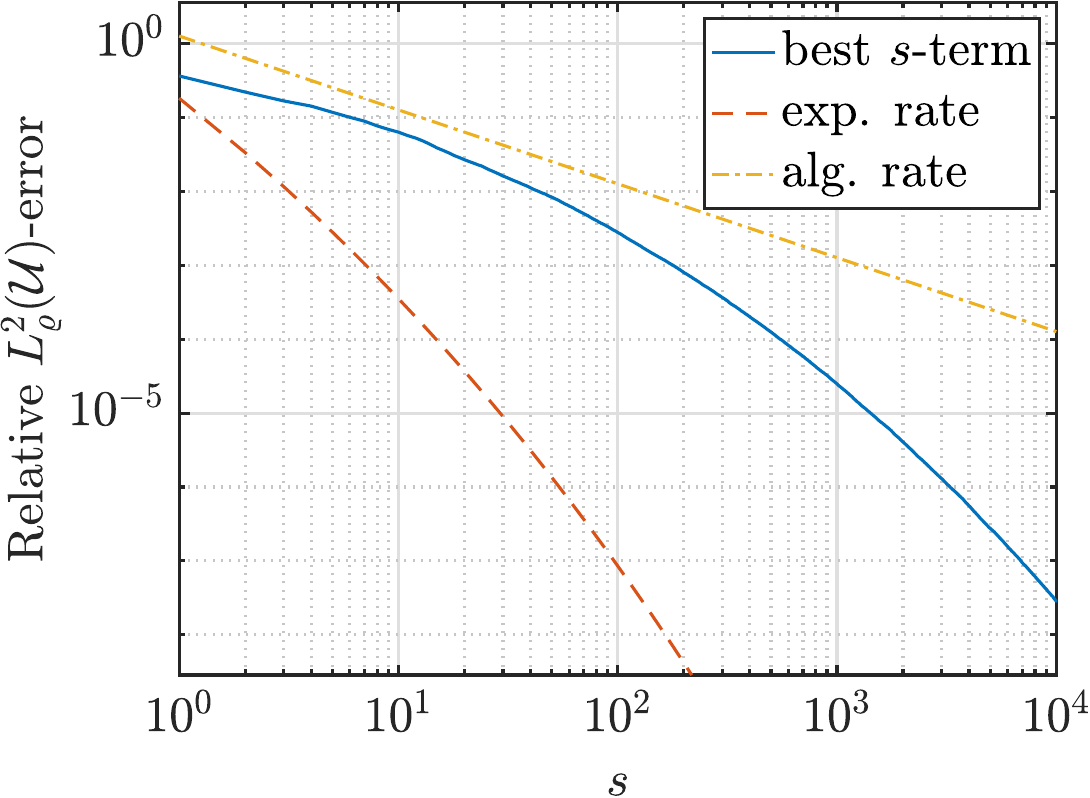}
\\
$d = 4$ & $d = 8$ 
\\
\includegraphics[width = 0.40\textwidth]{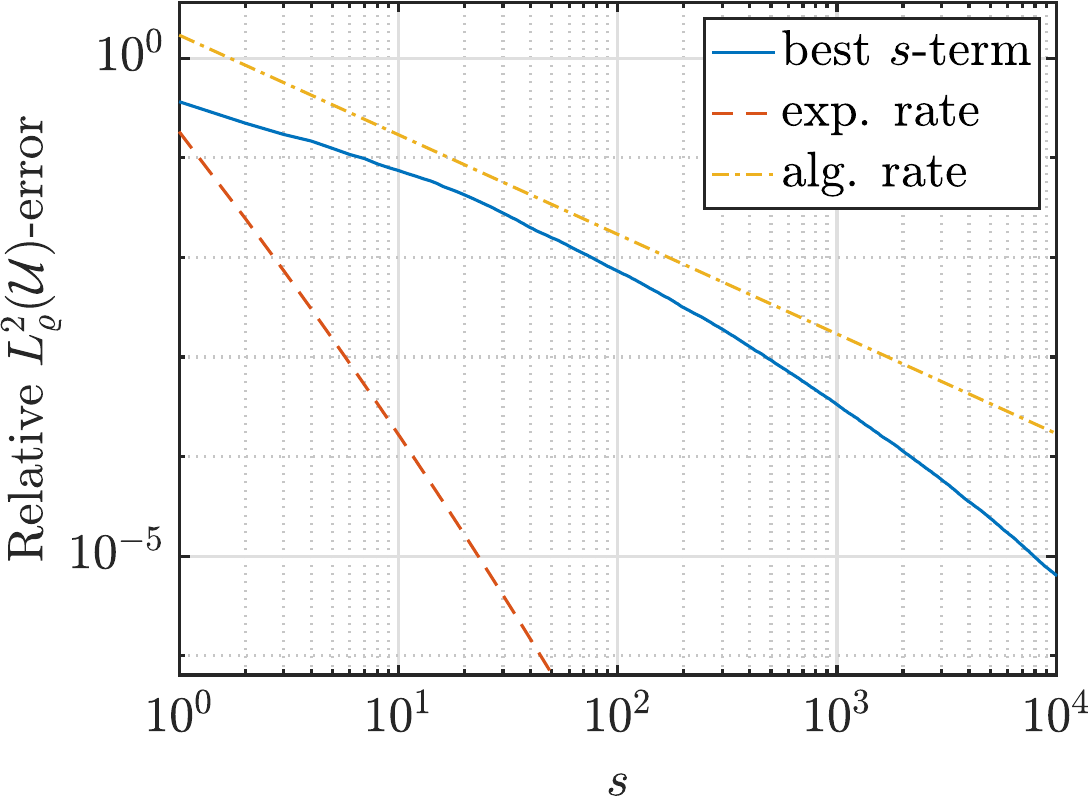}
&
\includegraphics[width = 0.40\textwidth]{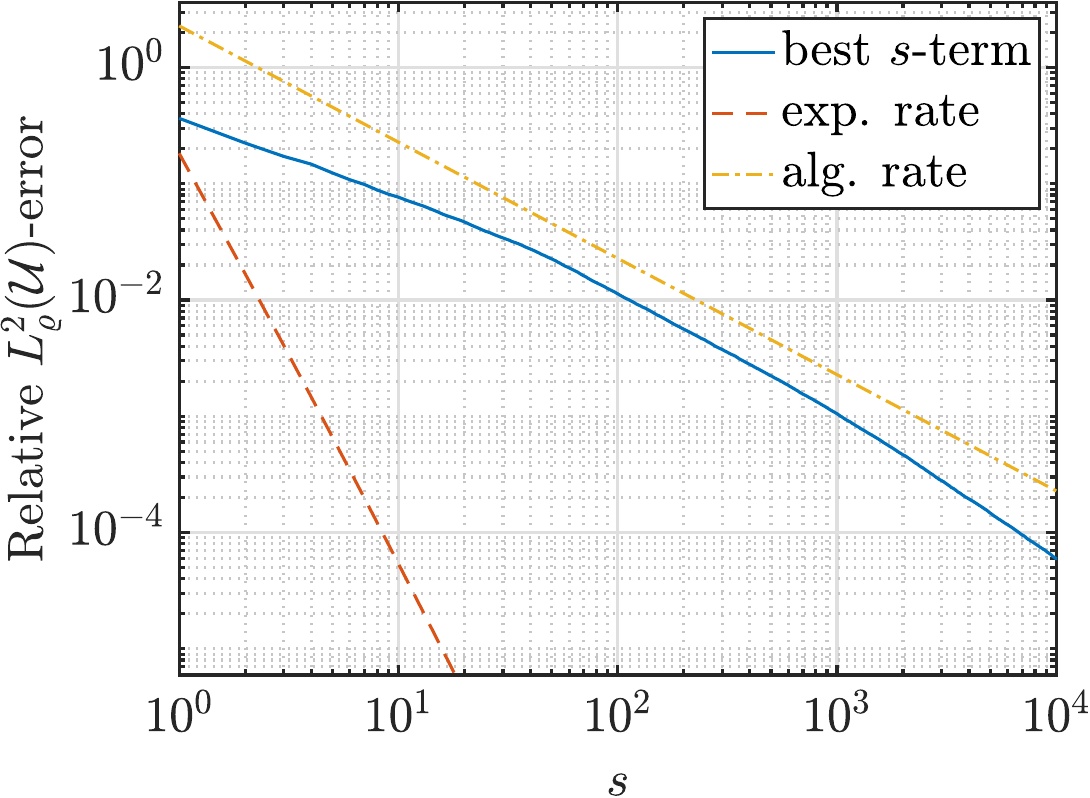}
\\
$d = 16$ & $d = 32$
\end{tabular}
\end{small}
\end{center}
\caption{Best $s$-term approximation error in the $L^2_{\varrho}(\cU)$-norm for \ef{f-numerics} with $\delta_i = i^{3/2}$. This figure also shows the exponential rate ``exp.\ rate'', defined as $C_{\mathsf{exp}} \cdot \exp \left ( - \left ( s d! \prod^{d}_{i=1} \log(\rho_i) \right )^{1/d} \right )$, where $\rho_i$ is such that $(\rho_i+1/\rho_i)/2 = 1+\delta_i$, and the algebraic rate ``alg.\ rate'', defined as  $C_{\mathsf{alg}} \cdot s^{-1}$. The constants $C_{\mathsf{exp}}$ and $C_{\mathsf{alg}}$ are chosen empirically to aid visualization.}
\label{fig:rates}
\end{figure}

 In this case, the algebraic rates of Theorem \ref{t:best-s-term} better describe the convergence behaviour. An example of this effect is shown in Figure \ref{fig:rates} for
 \be{
 \label{f-numerics}
 f(\bm{y}) = \prod^{d}_{i=1} (2 \delta_i + \delta^2_i)^{1/2} /  (y_i + 1 + \delta_i),\quad \forall \bm{y} = (y_i)^{d}_{i=1} \in [-1,1]^d.
 }
This function is holomorphic in $\cE(\bm{\rho})$ for any $\rho_i$ satisfying $(\rho_i + 1/\rho_i)/2  < 1+\delta_i$. Thus, by Remark \ref{rem:fin-dim}, it is also $(\bm{b},\varepsilon)$-holomorphic for any $b_i > \varepsilon/\delta_i$. In this figure, we consider $\delta_i = i^{3/2}$, meaning that $\bm{b}$ can be chosen so that $\bm{b} \in \ell^p(\bbN)$ for any $2/3 < p < 1$. Thus, Theorem \ref{t:best-s-term} predicts dimension-independent algebraic convergence with order that is arbitrarily close to $1/2-1/(2/3) = -1$. As we see, in $d = 4$ dimensions the error follows the exponential rate. However, once $d = 16$ or $d = 32$, the algebraic rate better predicts the true error.

This discussion and example motivates the study of algebraic rates. Simply put, they typically better describe the convergence behaviour, even in finite dimensions as long as $d$ is not too small (e.g., $d \geq 16$).

\section{Limits of learnability from data}\label{s:learn-limits}

We now consider learning $(\bm{b},\varepsilon)$-holomorphic functions from data. In this section, we present lower bounds for the amount of data that is necessary to learn such functions. We do this by appealing to techniques from information-based complexity \cite{novak1988deterministic,novak2008tractability,novak2010tractability,traub1988information} -- in particular, \textit{(adaptive) $m$-widths}.

Since our present focus is on lower bounds, in this section we do not assume that the training data takes the form of pointwise evaluations of the target function, as in \ef{training-data}. In fact, we can allow for arbitrary (adaptive) linear sampling operators. We will, for simplicity, consider only scalar-valued function approximation, i.e., $(\cV,\nm{\cdot}_{\cV}) = (\bbR,\abs{\cdot})$. See Remark \ref{rem:extensions}.

\subsection{Adaptive $m$-widths}

Let $C(\cU)$ be the Banach space of continuous functions $f : \cU \rightarrow \bbR$. We define an \textit{adaptive sampling operator} as a map of the form
\bes{
\cS : C(\cU)  \rightarrow \bbR^m,\ \quad \cS(f) = \begin{bmatrix} S_1(f) \\ S_2(f ; S_1(f)) \\ \vdots \\ S_m(f ; S_{1}(f),\ldots, S_{m-1}(f)) \end{bmatrix},
}
where  $S_1 : C(\cU) \rightarrow \bbR$ is a bounded  linear functional and, for $i = 2,\ldots,m$, $S_i :C(\cU) \times \bbR^{i-1} \rightarrow \bbR$ is bounded and linear in its first component.

This definition includes standard bounded linear operators as a special case. However, it also allows for situations where the $i$th sample is selected adaptively based on the existing measurements. In machine learning, this is commonly referred to as \textit{active learning} \cite{settles2012active}. In information-based complexity, this is commonly referred to as \textit{adaptive information} \cite[Sec.~4.1.1]{novak2008tractability}. As noted, our primary concern in this work is where each $S_i$ is a pointwise evaluation operator (so-called \textit{standard information} \cite[Sec.\ 4.1.1]{novak2008tractability}). In this case, $\cS(f) = (f(\bm{y}_i))^{m}_{i=1} \in \bbR^m$, where $\bm{y}_i \in \cU$ is the $i$th sample point, which is potentially chosen adaptively based on the previous measurements $f(\bm{y}_1),\ldots,f(\bm{y}_{i-1})$.

\defn{[Adaptive $m$-width]
\label{def:m-width}
The \textit{(adaptive) $m$-width} of a subset $\cK \subseteq C(\cU)$ is
\be{
\label{Theta_m-def}
\Theta_m(\cK) =  \inf_{\cS,\cR} \sup_{f \in \cK} \nm{f - \cR ( \cS(f)) }_{L^2_{\varrho}(\cU )},
}
where the infimum is taken over all adaptive sampling operators $\cS : C(\cU) \rightarrow \bbR^m$ and reconstruction maps $\cR : \bbR^m \rightarrow C(\cU)$. 
}

The adaptive $m$-width is related to the concept of \textit{information complexity}  \cite[\S4.1.4]{novak2008tractability}. It measures how well one can approximate functions from $\cK$ using an arbitrary (adaptive) linear sampling operator $\cS$ and an arbitrary (potentially nonlinear) reconstruction map $\cR$. Due to the inner supremum, the approximation is measured in a worst-case (or uniform) sense; in other words, sampling-recovery pair $(\cS,\cR)$ is required to provide a guaranteed bound simultaneously for all functions in the class.

\subsection{Lower bounds for adaptive $m$-widths}

We now present the main results of this section. Following \S \ref{ss:holo-p-classes}, we let
\bes{
\theta_m(p) = \Theta_m(\cH(p) ), \quad 
\theta_m(p,{\mathsf{M}}) = \Theta_m(\cH(p,{\mathsf{M}})).
}
 We also consider slightly weaker notions, which are defined as follows.
\bes{
\begin{split}
\overline{\theta_m}(p) &= \sup \left \{ \Theta_m(\cH(\bm{b})) : \bm{b} \in \ell^p(\bbN),\ \bm{b} \in [0,\infty)^{\bbN},\  \nm{\bm{b}}_p \leq 1 \right \},
\\
\overline{\theta_m}(p,{\mathsf{M}}) &= \sup \left \{ \Theta_m(\cH(\bm{b})) : \bm{b} \in \ell^p_{\mathsf{M}}(\bbN),\ \bm{b} \in [0,\infty)^{\bbN},\  \nm{\bm{b}}_{p,\mathsf{M}} \leq 1 \right \}.
\end{split}
}
The widths $\theta_m(p)$ and $\theta_m(p,{\mathsf{M}})$ pertain to the unknown anisotropy setting. A sampling-recovery pair $(\cS,\cR)$ does not have access to the anisotropy parameter $\bm{b}$, and must work uniformly for all holomorphic functions with $\ell^p$ or $\ell^p_{\mathsf{M}}$-summable $\bm{b}$ (of norm at most one). Conversely, $\overline{\theta_m}(p)$ and $\overline{\theta_m}(p,{\mathsf{M}})$ are weaker and pertain to the known anisotropy setting, since they allow the sampling-recovery pair $(\cS,\cR)$ to depend on $\bm{b}$. In particular, we have
\bes{
\theta_m(p) \geq \overline{\theta_m}(p)\text{ and }\theta_m(p,\mathsf{M}) \geq \overline{\theta_m}(p,\mathsf{M}).
}
\begin{theorem}[Limits of learnability]
\label{t:lower-bounds}
Let $m\geq 1$ and $0<p<1$. Then
\begin{itemize}
\item the $m$-widths  $\overline{\theta_m}(p)$ and $\overline{\theta_m}(p,\mathsf{M})$ satisfy
\begin{equation}\label{known-lower}
\overline{\theta_m}(p) \geq \overline{\theta_m}(p,\mathsf{M}) \geq c \cdot 2^{-\frac1p} \cdot m^{\frac12-\frac1p},
\end{equation}
\item the $m$-width $\theta_m(p)$ satisfies
\begin{equation} 
\label{unknown-lower-1}
\theta_m(p) \geq c \cdot 2^{\frac12-\frac2p},
\end{equation}
\item and the $m$-width $\theta_m(p,\mathsf{M})$ satisfies
\begin{equation}
\label{unknown-lower-2}
\theta_m(p,\mathsf{M}) \geq \overline{\theta_m}(p,\mathsf{M}) \geq c \cdot 2^{-\frac1p} \cdot m^{\frac12-\frac1p}.
\end{equation}
\end{itemize}
Here $c > 0$ is a universal constant.
\end{theorem}

See \cite[Thms.\ 4.4 \& 4.5]{adcock2023optimal}. This theorem has several consequences. First, \ef{unknown-lower-1} shows that it is impossible to uniformly learn functions from the class $\cH(p)$. In other words, mere $\ell^p$-summability of the anisotropy parameter $\bm{b}$ is insufficient. By contrast, \ef{unknown-lower-2} shows that it may be possible to learn functions in $\cH(p,\mathsf{M})$, where $\bm{b}$ is now $\ell^p_{\mathsf{M}}$-summable, but only at a rate that is algebraically decaying in $m$. We will see in the next section that these rates can in fact be attained (up to log terms) by practical methods. Finally, \ef{known-lower} asserts that knowledge of the anisotropy parameter $\bm{b}$ does not help. The lower bound for the nonuniform width $\overline{\theta_m}(p,\mathsf{M})$ is the same as the lower bound for the uniform width $\theta_m(p,\mathsf{M})$. Therefore, in terms of sample complexity, knowledge of $\bm{b}$ conveys no benefit -- we may as well deal exclusively with the unknown anisotropy setting.

It is worth relating this result back to Theorem \ref{t:best-s-term}. Theorem \ref{t:lower-bounds} shows that we can at best achieve a rate of $m^{1/2-1/p}$ when learning functions in $\cH(\bm{b})$ from data, whenever $\bm{b} \in \ell^p_{\mathsf{M}}(\bbN)$. Theorem \ref{t:best-s-term} (with $p = 2$) shows that the best $s$-term polynomial approximation can achieve a rate of $s^{1/2-1/p}$ subject to the weaker assumption $\bm{b} \in \ell^p(\bbN)$. Hence, the discrepancy between $\ell^p(\bbN)$ (where learning is impossible) and $\ell^p_{\mathsf{M}}(\bbN)$ (where learning is possible) must stem from the limited amount of data, not intrinsic properties of the function class itself.

\rem{[Extensions]
\label{rem:extensions}
For simplicity, we have limited this discussion to scalar-valued functions with $\varrho$ being the uniform measure. The latter assumption can be easily relaxed to an arbitrary tensor-product probability measure. One may also consider Banach-valued functions, although some care is needed when defining an adaptive sampling operator in this case.  See \cite{adcock2023optimal} for details.
}

\subsection{Towards methods}

We have now shown that functions in $\cH(p,\mathsf{M})$ can be learned from training sets of size $m$ with rates that are at best $\ord{m^{1/2-1/p}}$, regardless of the training data and learning method. In the remainder of this review, we focus on describing methods that achieve (close to) these rates when the training data consists of i.i.d.\ pointwise samples, as in \ef{training-data}. We consider, firstly, polynomial-based methods  (\S\ref{s:sparse-poly-learn}) and, secondly, methods based on DNNs (\S \ref{s:dnn-existence}-\ref{s:learning-proof}).

\section{Learning sparse polynomial approximations from data}\label{s:sparse-poly-learn}

We know from Theorem \ref{t:best-s-term} that the best $s$-term approximation converges at the desired rate $s^{1/2-1/p}$. Therefore, our goal is to design methods that can compute best (or quasi-best) $s$-term approximations using roughly $s$ samples (up to constant and log factors). We shall do this using tools from compressed sensing \cite{foucart2013mathematical,vidyasagar2019introduction}, which recast the problem of learning sparse polynomial approximations as recovering an approximately sparse vectors. Unfortunately, as we describe in this section, task is easily laid out, but not so easily completed. 

\subsection{Setup}\label{ss:CS-setup}

Throughout this section, we assume that $\cV$ is a Hilbert space (see \S \ref{ss:extensions} for some discussion on the Banach case). Moreover, 
in a practical scenario, we often cannot work directly with the space $\cV$, since it may be infinite dimensional. Therefore, we now consider a finite-dimensional discretization $\cV_h \subseteq \cV$,
where $h > 0$ denotes a discretization parameter. For parametric PDEs, this is typically the mesh size when a Finite Element Method (FEM) is used. We also assume that the training data \ef{training-data} belongs to $\cU \times \cV_h$, rather than $\cU \times \cV$, i.e.,
\be{
\label{training-data-Vh}
\{ (\bm{y}_i , f(\bm{y}_i) + e_i) \}^{m}_{i=1} \subseteq \cU \times \cV_h.
}
This encapsulates the notion that $f$ is evaluated using some numerical simulation (e.g., a FEM) that outputs values in $\cV_h$. Notice that the error term $e_i$ incorporates any error incurred by this computation. For convenience, we also define the orthogonal projection
\bes{
\cP_h : \cV \rightarrow \cV_h.
}
Moreover, if $f \in L^2_{\varrho}(\cU ; \cV)$, then we let $\cP_h f \in L^2_{\varrho}(\cU ; \cV_h)$ be the function defined almost everywhere as $(\cP_h f)(\bm{y}) = \cP_h(f(\bm{y}))$, $\bm{y} \in \cU$.

Now consider a function $f \in L^2_{\varrho}(\cU ; \cV)$ with expansion \ef{f_exp} and coefficients \ef{f-coeff}. Standard compressed sensing involves the recovery of finite vectors. Therefore, we first need to truncate the expansion \ef{f_exp}. Let $\Lambda \subset \cF$ be finite with $|\Lambda|=N$ and an enumeration $\Lambda = \{\bm{\nu}_1,\ldots,\bm{\nu}_N \}$. Then
\bes{
f(\bm{y}_i) + e_i = f_{\Lambda}(\bm{y}_i) + (f-f_{\Lambda})(\bm{y}_i) +e_i = \sum_{\bm{\nu} \in \Lambda} c_{\bm{\nu}} \Psi_{\bm{\nu}}(\bm{y}_i) + (f-f_{\Lambda})(\bm{y}_i) + e_i,
}
where $f_{\Lambda} = \sum_{\bm{\nu} \in \Lambda} c_{\bm{\nu}} \Psi_{\bm{\nu}}$. Now let  $\bm{c}_{\Lambda} = (c_{\bm{\nu}_i})^{N}_{i=1} \in \cV^N$. Then we have
\be{
\label{AcLambda}
\bm{A} \bm{c}_{\Lambda} + \bm{e} + \bm{e}' = \bm{f},
}
where 
\be{
\label{Ab-CS-def}
\bm{f} = \frac{1}{\sqrt{m}} \left ( f(\bm{y}_i) + e_i  \right )^{m}_{i=1}
\quad 
\bm{A} = \frac{1}{\sqrt{m}} \left (  \Psi_{\bm{\nu}_j}(\bm{y}_i)/\sqrt{m} \right )^{m,N}_{i,j=1} \in \bbR^{m \times N},
}
and
\bes{
\bm{e} = \frac{1}{\sqrt{m}} (e_i)^{m}_{i=1},\quad \bm{e}' = \frac{1}{\sqrt{m}} \left ( (f-f_{\Lambda})(\bm{y}_i)  \right )^{m}_{i=1} \in \bbR^m.
}
The vector $\bm{c}_{\Lambda}$ is approximately $s$-sparse. Therefore, we have recast the problem into a standard compressed sensing form: namely, the recovery of an approximately sparse vector from noisy linear measurements \ef{AcLambda}.

Of course, this problem is not quite standard, since the vector $\bm{c}_{\Lambda}$ has entries taking values in $\cV$, rather than $\bbR$ or $\bbC$. Fortunately, this challenge can be dealt with by `lifting' the appropriate theoretical tools from $\bbR$ or $\bbC$ to Hilbert spaces (see Remark \ref{rem:lifting} below). A more delicate challenge is the following.

\pbk
\textbf{Challenge 1.}\ How should the truncation set $\Lambda$ be chosen.

\pbk
In general, $f$ has an infinite expansion. Since we recover only those coefficients with indices in $\Lambda$ (i.e., the vector $\bm{c}_{\Lambda}$), the error in recovering $f$ will always involve the expansion tail $f - f_{\Lambda}$. On the one hand, $\Lambda$ should not be too large, since, as we see later in \ef{samp-comp-CS}, $\log(N)$ will enter into the sample complexity bound. Moreover, larger $N$ will also increase the computational cost. On the other hand, we need $f - f_{\Lambda}$ to be sufficiently small so as to obtain optimal learning rates. To do this, we need to ensure no large coefficients of $f$ lie outside $\Lambda$. And here is where the problem lies. The result on best $s$-term approximation, Theorem \ref{t:best-s-term}, which has served as our rationale up to now for using compressed sensing, asserts approximate sparsity of the coefficients, but gives not guarantees on where the largest $s$ coefficients should lie within the infinite set $\cF$.

\subsection{Sampling discretizations for multivariate polynomials}\label{ss:sampling-disc}

We will return to Challenge 1 in a moment. But, first, we also need to introduce a second challenge.
Let us imagine some oracle gives us a `good' set $S \subset \cF$ of size $|S| = s$. For example, this could even be a set $S = S^*$ corresponding to the best $s$-term approximation \ef{best-s-term-inf}. Then one would naturally learn an approximation to $f$ via the empirical least-squares fit
\be{
\label{oracle-LS}
\hat{f} \in \argmin{p \in \bbP_{S;\cV}} \frac1m \sum^{m}_{i=1} \nm{ f(\bm{y}_i) + e_i - p(\bm{y}_i) }^2_{\cV},
}
where $\bbP_{S;\cV} $ is the $s$-dimensional subspace
\bes{
\bbP_{S;\cV} = \left \{ \sum_{\bm{\nu} \in S} c_{\bm{\nu}} \Psi_{\bm{\nu}} : c_{\bm{\nu}} \in \cV \right \} \subset L^2_{\varrho}(\cU ; \cV)
}
of Hilbert-valued polynomials with non-zero coefficients in $S$. 

The behaviour of the estimator \ef{oracle-LS} is intimately related to the existence of a \textit{sampling discretization} of the $L^2_{\rho}(\cU)$-norm for the scalar-valued analogue of this space, i.e., $\bbP_{S} = \bbP_{S;\bbR} = \spn \{ \Psi_{\bm{\nu}} : \bm{\nu} \in \cS \}$.
A sampling discretization  \cite{kashin2022sampling} (also known as a \textit{Marcinkiewicz--Zygmund inequality} \cite{temlyakov2018marcinkiewicz}) is an inequality of the form
\be{
\label{samp-disc}
\alpha \nm{p}^2_{L^2_{\varrho}(\cU)} \leq \frac1m \sum^{m}_{i=1} | p (\bm{y}_i) |^2 \leq \beta \nm{p}^2_{L^2_{\varrho}(\cU)},\quad \forall p \in \bbP_{S},
}
for constants $0 < \alpha \leq \beta < \infty$. Specifically, \ef{samp-disc} ensures both accuracy of the estimator $\hat{f}$ and robustness to measurement error (see, e.g., \cite[Thm.\ 5.3]{adcock2022sparse}).

\rem{[Lifting to Hilbert spaces]
\label{rem:lifting}
This is an instance of the aforementioned `lifting' concept. The sampling discretization \ef{samp-disc} is formulated for the space $\bbP_{S} = \bbP_{S;\bbR}$. It is a short argument (see, e.g., \cite[Lem.\ 7.5]{adcock2024efficient}) to show that it is equivalent to a sampling discretization for the space $\bbP_{S;\cV}$, i.e.,
\bes{
\alpha \nm{p}^2_{L^2_{\varrho}(\cU;\cV)} \leq \frac1m \sum^{m}_{i=1} \nm{p (\bm{y}_i) }^2_{\cV} \leq \beta \nm{p}^2_{L^2_{\varrho}(\cU;\cV)},\quad \forall p \in \bbP_{S;\cV}.
}
In turn, this implies accuracy and robustness of the estimator \ef{oracle-LS}.
}

Sufficient conditions for sampling discretizations in linear subspaces with i.i.d.\ samples can be derived using standard matrix concentration inequalities. See, e.g., \cite{cohen2013stability} or \cite[Chpt.\ 5]{adcock2022sparse}. These conditions involve the quantity
\be{
\label{kappa-PS}
\kappa(\bbP_S) = \nm{\cK(\bbP_S)}_{L^{\infty}_{\rho}(\cU)},
}
where $\cK(\bbP_S)$ is the (reciprocal) \textit{Christoffel function} of the subspace $\bbP_S$:
\be{
\label{K-PS}
\cK(\bbP_S)(\bm{y}) = \sum_{\bm{\nu} \in S} | \Psi_{\bm{\nu}}(\bm{y}) |^2,\quad \forall \bm{y} \in \cU.
}
Specifically, one can show that if
\be{
\label{samp-comp-LS}
m \geq c \cdot \kappa(\bbP_S) \cdot \log(2s/\epsilon),
}
where $s = |S|$ and $c > 0$ is a universal constant, then \ef{samp-disc} holds with constants $\beta \leq 2$ and $\alpha \geq 1/2$ (these values are arbitrary) with probability at least $1-\epsilon$ on the draw of the sample points $\bm{y}_1,\ldots,\bm{y}_m \sim_{\mathrm{i.i.d.}} \varrho$ (see, e.g., \cite[Thm.\ 5.12]{adcock2022sparse}).

Therefore, the \textit{sample complexity} of the `oracle' estimator \ef{oracle-LS} is governed by the maximal behaviour of the Christoffel function \ef{kappa-PS}. Note that $\kappa(\bbP_S) \geq s$ for any set $S$ with $|S| = s$ \cite[\S 5.3]{adcock2022sparse}. To use \ef{samp-comp-LS} to achieve the optimal rates described in Theorem \ref{t:lower-bounds} for this oracle estimator, we require that $\kappa(\bbP_S) \lesssim s$ as well. Unfortunately, this is not the case: $\kappa(\bbP_S)$ can be arbitrarily large in comparison to $s$. To see why, we recall \ef{P-max}. This implies that
\be{
\label{kappaPS-char}
\kappa(\bbP_S) = \cK(P_{S})(\bm{1}) = \sum_{\bm{\nu} \in S} u^2_{\bm{\nu}},
}
where the $u_{\bm{\nu}}$ are as in \ef{u-def}. Since $u_{\bm{\nu}} \rightarrow \infty$ as $\bm{\nu} \rightarrow \infty$, we deduce the claim.

\pbk
\textbf{Challenge 2.} Even if it is known, the index set $S$ of the best $s$-term approximation may lead to a sample complexity estimate \ef{samp-comp-LS} that is arbitrarily large.

\pbk
Similar to Challenge 1, this difficulty arises because Theorem \ref{t:best-s-term} gives no guarantees about the set $S = S^*$ which attains the best $s$-term approximation. In particular, it says nothing about the term $\kappa(\bbP_{S^*})$. The paths to resolving each challenge are therefore similar. We need to show that \textit{near-best} $s$-term approximations can be obtained using suitably \textit{structured} index sets.

\rem{
[Necessity of sampling discretizations]
The sampling discretization \ef{samp-disc} is sufficient condition for robustness of the estimator $\hat{f}$. However, the lower inequality is also essentially necessary. Indeed, let $\cS : C(\cU) \rightarrow \bbR^m, f \mapsto (f(\bm{y}_i))^{m}_{i=1} / \sqrt{m}$ and $\cR : \bbR^m \rightarrow C(\cU)$ be any reconstruction map that is \textit{$\delta$-accurate} over $\bbP_S$, i.e., 
\be{
\label{acc-rec-map}
\nm{p - \cR(\cS(p))}_{L^2_{\varrho}(\cU)} \leq \delta \nm{p}_{L^2_{\varrho}(\cU)},\quad \forall p \in \bbP_s,
}
for some $\delta > 0$.
Then it is a short argument to show that the $\epsilon$-Lipschitz constant
\bes{
L_{\epsilon} = \sup_{f \in C(\cU)} \sup_{\substack{0 < \nm{e}_2 \leq \epsilon}}  \nm{\cR(\cS(f)+e) - \cR(\cS(f))}_{L^2_{\varrho}(\cU)} / \nm{e}_2
}
satisfies $L_{\epsilon} \geq (1-\delta)\sqrt{\alpha}$. Thus, when $\alpha$ is small, a reconstruction map cannot be simultaneously accurate (in the sense of \ef{acc-rec-map}) and stable.
}

\subsection{Resolving Challenge 1: lower and anchored sets}\label{ss:challenge1-resolve}

In order to address Challenge 1, we now introduce the following concept.

\defn{[Lower and anchored sets]
A multi-index set $S \subseteq \cF$ is \textit{lower} if the following holds for every $\bm{\nu}, \bm{\mu} \in S$:
\bes{
(\bm{\nu} \in S\text{ and } \bm{\mu} \leq \bm{\nu}) \Rightarrow \bm{\mu} \in \Lambda.
}
A multi-index set $S \subseteq \cF$ is \textit{anchored} if it is lower and if the following holds for every $j \in \bbN$:
\bes{
\bm{e}_j \in S \Rightarrow \{\bm{e}_1,\bm{e}_2, \ldots, \bm{e}_{j} \}\subseteq S.
}
}
Lower sets are classical objects in multivariate approximation theory. Anchored sets were introduced in the context of infinite-dimensional approximations. See, e.g., \cite[\S 2.3.3 \& 2.5.3]{adcock2022sparse} and references therein. 

We now state a result that shows that the $s^{1/q-1/p}$ rate asserted in Theorem \ref{t:best-s-term} can also be obtained using a lower or (under an additional assumption) anchored set. See, e.g., \cite[Thm.\ 3.33]{adcock2022sparse} or \cite[\S 3.8]{cohen2015approximation}.

\thm{
[Algebraic convergence in lower or anchored sets]
\label{t:best-s-term-anchored}
Let $\bm{b} \in [0,\infty)^{\bbN}$ be such that $\bm{b} \in \ell^p(\bbN)$ for some $0 < p < 1$. Then for any $s \in \bbN$ and $p \leq q \leq 2$, there exists an lower set $S \subset \cF$ with $|S|\leq s$ such that
\bes{
\nm{\bm{c} - \bm{c}_S}_{q;\cV} \leq C(\bm{b},p) \cdot s^{\frac1q-\frac1p}
}
for all $f \in \cH(\bm{b})$ with coefficients $\bm{c}$ as in \ef{f-coeff}. If $\bm{b} \in \ell^p_{\mathsf{M}}(\bbN)$ then $S$ can also be chosen as an anchored set.
}

In particular, when $q = 2$, this theorem and Parseval's identity imply the existence of a lower or anchored set $S$ with $|S| \leq s$ such that
\be{
\label{anchored-rate}
\nm{f - f_S}_{L^{2}_{\varrho}(\cU ; \cV)} \leq C(\bm{b},p) \cdot s^{\frac12-\frac1p},\quad \forall f \in \cH(\bm{b}).
}
Now recall that our aim is to learn functions in the space $\cH(p,\mathsf{M})$ defined by \ef{HpM}. Any $f \in \cH(p,\mathsf{M})$ satisfies $f \in \cH(\bm{b})$ for some $\bm{b} \in \ell^p_{\mathsf{M}}(\bbN)$. Therefore, for any such function, we know there is an anchored set $S$, $|S| \leq s$, that achieves \ef{anchored-rate}.  How does this allow us to overcome Challenge 1? The reason is because anchored sets of a fixed size lie within a \textit{finite} subset of $\cF$. Indeed, one can show (see, e.g., \cite{cohen2017discrete} or \cite[Prop.\ 2.18]{adcock2022sparse}) that $S \subset \Lambda^{\mathsf{HCI}}_{s}$ for all $S \subset \cF$ with $|S| \leq s$, where $\Lambda^{\mathsf{HCI}}_{s}$ is the finite set
\be{
\label{HCI-def}
\Lambda^{\mathsf{HCI}}_{s} = \left \{ \bm{\nu} = (\nu_k)^{\infty}_{k=1} \in \cF : \prod^{s-1}_{k=1} (\nu_k+1) \leq s,\ \nu_k = 0,\ \forall k \geq s \right \}.
}
This set is in fact isomorphic to the $(s-1)$-dimensional \textit{hyperbolic cross} set of order $s-1$, a very well-known object in high-dimensional approximation \cite{dung2018hyperbolic}. With this knowledge, we set $\Lambda = \Lambda^{\mathsf{HCI}}_s$ and then apply Theorem \ref{t:best-s-term-anchored}:
\bes{
\nm{f - f_{\Lambda}}_{L^2_{\varrho}(\cU ; \cV)} = \nm{\bm{c} - \bm{c}_{\Lambda}}_{2;\cV} \leq \nm{\bm{c} - \bm{c}_S}_{2;\cV} \leq C(\bm{b},p) \cdot s^{\frac12-\frac1p},
}
for all $f \in \cH(\bm{b})$ and $\bm{b} \in \ell^p_{\mathsf{M}}(\bbN)$. This resolves Challenge 1.

\rem{
Theorem \ref{t:lower-bounds} states that no method can learn functions in $\cH(p)$ from finite data. Theorem \ref{t:best-s-term-anchored} says that we can always find a lower set that yields the desired rate \ef{anchored-rate} when $\bm{b} \in \ell^p(\bbN)$. However, the union of all lower sets in infinite dimensions is not a finite set. This is precisely what prohibits one from learning functions in $\cH(p)$ using compressed sensing (which, of course, must not be possible in view of Theorem \ref{t:lower-bounds}): namely, there is no way to construct a suitable $\Lambda$ that ensures a uniformly small truncation error $f - f_{\Lambda}$.
}

\subsection{Resolving Challenge 2: weighted $k$-term approximation}

We now return to Challenge 2. In the previous subsection, we identified anchored sets as structured sets which achieve near-best $s$-term approximation rates. Unfortunately, while anchored sets do ameliorate the sample complexity issue, they do not fully resolve it. On the one hand, it is possible to show that  $\kappa(\bbP_S) \leq s^2$ whenever $S$ is anchored (see, e.g., \cite[Prop.\ 5.17]{adcock2022sparse}), where  $\kappa(\bbP_S)$ is as in \ef{kappa-PS}. Unfortunately, this bound is sharp. For example, the anchored set 
$$
S = \{ j \bm{e}_1 : j = 0,\ldots,s-1\},\quad \text{where }\bm{e}_1 = (1,0,0,\ldots),
$$
satisfies $\kappa(\bbP_S) = s^2$. This follows from \ef{kappaPS-char} and the fact that $u_{j \bm{e}_1} = u_j = \sqrt{2j+1}$. Therefore, the sample complexity bound \ef{samp-comp-LS} for learning with anchored sets is generically log-quadratic in $s$. 
 
To overcome Challenge 2 we need a different concept, \textit{weighted $k$-term approximation}. Motivated by \ef{kappaPS-char}, we now define the \textit{weighted cardinality} of set $S \subset \cF$
with respect to \textit{weights} $\bm{u}$ as
\bes{
| S |_{\bm{u}} : = \sum_{\bm{\nu} \in S} u^2_{\bm{\nu}} .
}
Thus we may reinterpret \ef{samp-comp-LS} as follows: the sample complexity is determined not by the cardinality of $S$, but by its weighted cardinality. 

Fortunately, as we next show, there exist sets of a given weighted cardinality that achieve near-optimal approximation rates. For this, we need some additional notation. Given $1 \leq p \leq 2$, we define the weighted $\ell^p_{\bm{u}}({\Lambda}; \cV)$ space as the space of $\cV$-valued sequences $\bm{c} = (c_{\bm{\nu}})_{\bm{\nu} \in {\Lambda}}$ for which $\nm{\bm{c}}_{p,\bm{u};\cV} < \infty$, where 
\bes{
\nm{\bm{c}}_{p,\bm{u};\cV} =  \left ( \sum_{\bm{\nu} \in {\Lambda}} u^{2-p}_{\bm{\nu}} \nm{c_{\bm{\nu}}}^p_{\cV} \right )^{\frac1p} .
}
Notice that $\| \cdot \|_{2,\bm{u} ; \cV}$ coincides with the unweighted norm $\| \cdot \|_{2;\cV}$.

\thm{
[Algebraic convergence of the weighted best $k$-term approximation]
\label{t:best-s-term-weighted}
Let $\bm{b} \in [0,\infty)^{\bbN}$ be such that $\bm{b} \in \ell^p(\bbN)$ for some $0 < p < 1$. Then for any $k > 0$ and $p \leq q \leq 2$ there exists a set $S \subset \cF$ with $|S|_{\bm{u}} \leq k$ such that
\bes{
\nm{\bm{c} - \bm{c}_S}_{q,\bm{u};\cV} \leq C(\bm{b},p) \cdot k^{\frac1q-\frac1p}
}
for all $f \in \cH(\bm{b})$ with coefficients $\bm{c}$ as in \ef{f-coeff}.
}

This theorem shows that we can construct near-best approximations using set of weighted cardinality at most $k$. In particular, when $q = 2$, this theorem and Parseval's identity give that there exists a set $S$ with $|S|_{\bm{u}} \leq k$ such that
\bes{
\nm{f - f_{S}}_{L^{2}_{\varrho}(\cU ; \cV)} \leq C(\bm{b},p) \cdot k^{\frac12-\frac1p},\quad \forall f \in \cH(\bm{b}).
}
Note that we do not require $\bm{b} \in \ell^p_{\mathsf{M}}(\bbN)$ for this result. This additional regularity is only needed to resolve Challenge 1.

\subsection{Weighted $\ell^1$-minimization}

Challenge 2 is now resolved, at least insofar as the `oracle' least-squares estimator goes. Of course, in practice we do not have access to such an oracle. However, Theorem \ref{t:best-s-term-weighted} implies that the sequence $\bm{c}$, and therefore the finite vector $\bm{c}_{\Lambda}$, is approximately \textit{weighted sparse}. In other words, it is well approximated by $\bm{c}_{S}$ for some index set $S \subseteq \Lambda$ of weighted cardinality $|S|_{\bm{u}} = k$. 

Weighted sparsity is a well-understood concept in compressed sensing (see \cite{adcock2017infinite-dimensional,rauhut2016interpolation} and references therein). Much like how sparsity can be promoted by solving a minimization problem involving the $\ell^1$-norm, weighted sparsity can be solving a minimization problem involving the weighted $\ell^1_{\bm{u}}$-norm. A large weight penalizes the corresponding coefficient, much like how a large weight increases the weighted sparsity of any vector that has a nonzero coefficient at the corresponding index. There are a number of ways to formulated the weighted $\ell^1$-minimization problem, but, following \cite{adcock2019correcting,adcock2024efficient}, we will consider the Hilbert-valued, \textit{weighted square-root LASSO} program
\be{
\label{fdef1}
\min_{\bm{z} \in \cV_h^{N}} \cG(\bm{z}),\quad \text{where }\cG(\bm{z}) : = \lambda \nm{\bm{z}}_{1,\bm{u};\cV} + \nm{\bm{A} \bm{z} - \bm{f}}_{2;\cV}.
}
Here $\nm{\bm{z}}_{1,\bm{u};\cV} = \sum_{\bm{\nu} \in \Lambda} u_{\bm{\nu}} \nm{z_{\bm{\nu}}}_{\cV} $ is the $\ell^1_{\bm{u};\cV}$-norm and $\lambda > 0$ is a parameter. Notice that this problem is posed over $\cV^N_h$, which means it can be numerically solved. See Remark \ref{basis-practical-soln} below.

In order to account for inexact solution of \ef{fdef1}, given $\gamma \geq 0$ we say that $\hat{\bm{c}} = (\hat{c}_{\bm{\nu}})_{\bm{\nu} \in \Lambda}$ is a \textit{$\gamma$-minimizer} of $\ef{fdef1}$ if
\bes{
\cG(\hat{\bm{c}}) \leq \min_{\bm{z} \in \cV_h^{N}} \cG(\bm{z}) + \gamma.
}
For such $\hat{\bm{c}}$, we define the corresponding sparse polynomial approximation
\be{
\label{fdef2}
\hat{f} = \sum_{\bm{\nu} \in \Lambda} \hat{c}_{\bm{\nu}} \Psi_{\bm{\nu}}.
}

\rem{
\label{basis-practical-soln}
In practice, \ef{fdef1} is solved by first introducing a basis $\{ \varphi_i \}^{K}_{i=1}$ for the space $\cV_h$, where $K = \dim(\cV_h)$ (in particular, $K = 1$ in the scalar-valued case $\cV = (\bbR,\abs{\cdot})$) . Instead of searching for a $\cV_h$-valued vector of coefficients $\hat{\bm{c}} \in \cV^N_h$, one now searches for an equivalent matrix of coefficients $\widehat{\bm{C}} \in \bbR^{N \times K}$. It is a short exercise to show that \ef{fdef1} is equivalent to theproblem
\be{
\label{min-equiv}
\min_{\bm{Z} \in \bbR^{N \times K}} \lambda \nm{\bm{Z}}_{2,1,\bm{u}} + \nm{(\bm{A} \bm{Z} - \bm{F}) \bm{G}^{1/2}}_{F},
}
where $\nm{\cdot}_{F}$ is the Frobenius norm, $\bm{G} = \left ( \ip{\varphi_i}{\varphi_j} \right )^{K}_{i,j=1} \in \bbR^{K \times K}$and
\bes{
 \nm{\bm{Z}}_{2,1,\bm{u}} = \sum^{N}_{i=1} u_{\bm{\nu}_i} \sqrt{ \sum^{K}_{j=1} |z_{ij}|^2 }.
}
As discussed in \cite{adcock2024efficient}, the convex optimization problem \ef{min-equiv} can be solved efficiently using Chambolle \& Pock's primal-dual iteration \cite{chambolle2011first-order,chambolle2016ergodic} in combination with a restart scheme \cite{adcock2023restarts,roulet2020sharpness}. Specifically, one can obtain a $\gamma$-minimizer in at most $\log(1/\gamma)$ iterations, where the cost-per-iteration is bounded by
\be{
\label{comp-cost}
c \cdot \left ( m \cdot N \cdot K + (m+N) \cdot (F(\bm{G}) + K) \right ).
}
Here $F(\bm{G}) \leq K^2$ is the cost of performing the matrix-vector multiplication $\bm{x} \mapsto \bm{G} \bm{x}$. See \cite[Thm.\ 3.9 \& Lem.\ 4.3]{adcock2024efficient}.
}

\subsection{Theoretical guarantee}

We are now ready to present a theoretical guarantee for this estimator. For convenience, we define
\be{
\label{Ldef}
L = L(m,\epsilon) = \log( m) \cdot ( \log^3( m) + \log(\epsilon^{-1})  ).
}

\thm{
[Near-optimal learning via polynomials]
\label{thm:main-res-poly}
Let $m \geq 3$, $0 < \epsilon < 1$ and $n = \lceil m / L \rceil$, where $L = L(m,\epsilon)$ is as in \ef{Ldef}. Let $0 < p < 1$, $\bm{b} \in \ell^p_{\mathsf{M}}(\bbN)$, $f \in \cH(\bm{b})$, $\bm{y}_1,\ldots,\bm{y}_m \sim_{\mathrm{i.i.d.}}\varrho$ and consider the training data \ef{training-data-Vh}. Then with probability at least $1-\epsilon$ the approximation $\hat{f}$ defined by \ef{fdef2} for any $\gamma$-minimizer $\hat{\bm{c}}$ of \ef{fdef1}, $\gamma \geq 0$, with $\Lambda = \Lambda^{\mathsf{HCI}}_n$ satisfies
\ea{
\label{err-bound-map-alg-inf_1}
\nm{f - \hat{f}}_{L^2_{\varrho}(\cU;\cV)} \leq c \cdot \zeta,
}
where $c \geq 1$ is a universal constant,
\be{
\label{zeta-alg-inf-def}
\zeta := C \cdot \left ( m/L\right )^{1/2-1/p}  + \nm{\bm{e}}_{2;\cV} / \sqrt{m} +  \nm{f - \cP_h(f)}_{L^{\infty}(\cU ; \cV)} + \gamma.
}
}

See \cite[Thm.\ 3.7]{adcock2024efficient}. The above theorem is slightly more general than \cite[Thm.\ 3.7]{adcock2024efficient} since it does not require $\bm{b}$ to be monotonically decreasing. This generalization follows from techniques given in the proof of Theorem 7.2 in \cite{adcock2023monte}.

This result demonstrates that functions in $\cH(p,\mathsf{M})$ can be learned from i.i.d.\ samples with a rate $\ord{(m/\log^4(m))^{1/2-1/p}}$ which, in view of Theorem \ref{t:lower-bounds}, is near-optimal. Moreover, the estimator is robust to other sources of error in the problem.  The term $\zeta$ defined in \ef{zeta-alg-inf-def} divides into the following terms:
\begin{enumerate}
\item[(a)] An \textit{approximation error} $C \cdot (m/L)^{1/2-1p}$, as discussed.
\item[(b)] A \textit{measurement error} $\nm{\bm{e}}_{2;\cV}/\sqrt{m}$, which accounts for any errors in computing the sample values $f(\bm{y}_i)$.
\item[(c)] A \textit{physical discretization error} $\nm{f - \cP_h(f)}_{L^{\infty}(\cU ; \cV)}$, which accounts for the error induced when working over the finite-dimensional space $\cV_h$ instead of $\cV$ and depends on the orthogonal projection $\cP_h(f)$.
\item[(d)] An \textit{optimization error} $\gamma$, which depends on the optimality gap of the computed solution $\hat{\bm{c}}$. See Remark \ref{basis-practical-soln} for further discussion on this term.
\end{enumerate}
The proof of this theorem relies heavily on tools from compressed sensing. To relate it back to the discussion in \S \ref{ss:sampling-disc}, we observe that a key step in the proof is establishing a sampling discretization of the form
\be{
\label{universal-samp-disc}
\alpha \nm{p}^2_{L^2_{\varrho}(\cU)} \leq \frac1m \sum^{m}_{i=1} | p (\bm{y}_i) |^2 \leq \beta \nm{p}^2_{L^2_{\varrho}(\cU)},\quad \forall p \in \bbP_{S},\ S \subseteq \Lambda,\ |S|_{\bm{u}} \leq k.
}
This is stronger than \ef{samp-disc}, since it is required to hold simultaneously for all sets $S \subseteq \Lambda$ of weighted cardinality at most $k$. It is an example of a \textit{universal sampling discretization} \cite{dai2023universal}. The proof first shows that it is equivalent to a certain \textit{weighted Restricted Isometry Property (wRIP)}, then uses known results for the wRIP (see \cite{brugiapaglia2021sparse}, as well as \cite{rauhut2016interpolation}). As shown therein, a sufficient condition for \ef{universal-samp-disc} to hold with probability at least $1-\epsilon$ and constants $\alpha \geq 1-\delta$ and $\beta \leq 1+\delta$ for some $0 < \delta < 1$ is
\be{
\label{samp-comp-CS}
m \geq c \cdot \delta^{-2} \cdot k \cdot ( \log(\E N) \cdot \log^2(k/\delta) + \log(2/\epsilon) ).
}
Crucially, up to the log terms, there is now a linear scaling between $m$ and $k$. This is what leads to the near-optimal approximation error term.

\rem{
\label{rem:postprocess-sparse}
Generally, the vector $\hat{\bm{c}}$ computed by solving \ef{fdef1} will not be sparse. Therefore, $\hat{f}$ will not be a \textit{sparse polynomial} approximation per se. Fortunately, one can always postprocess the solution to obtain such an approximation. This was shown in \cite[Thm.\ 7.2]{adcock2023monte} in the scalar-valued case $\cV = (\bbR,\abs{\cdot})$, but the technique extends straightforwardly to the Hilbert-valued case. In the case of Theorem \ref{thm:main-res-poly}, given a $\gamma$-minimizer $\hat{\bm{c}} = (\hat{c}_{\bm{\nu}} )_{\bm{\nu} \in \Lambda} \in \cV^N_h$ one first computes the index set $S \subseteq \Lambda$ of the largest $n = \lceil m / L \rceil$ entries of the vector $(\nm{\hat{c}_{\bm{\nu}}}_{\cV})_{\bm{\nu}\in \Lambda}$, and then replaces $\hat{f}$ with the $n$-sparse polynomial approximation $\check{f} = \sum_{\bm{\nu} \in S} \hat{c}_{\bm{\nu}} \Psi_{\bm{\nu}}$. As shown in \cite[Thm.\ 7.2]{adcock2023monte}, $\check{f}$ attains the same error bounds up to possible changes in the constants. The computational cost of this postprocessing step, roughly $\ord{N \log(N) + N \cdot F(\bm{G})}$ operations, is generally negligible in comparison to the cost of computing the initial $\gamma$-minimizer $\hat{\bm{c}}$.
}

\subsection{Discussion and extensions}\label{ss:extensions}

We conclude this section a short discussion. First, Theorem \ref{thm:main-res-poly} asserts that i.i.d.\ pointwise samples constitute \textit{near-optimal information}. This is an interesting facet of infinite-dimensional holomorphic function approximation, which contrasts starkly with the finite-, and specifically, low-dimensional case, where samples drawn i.i.d.\ from the uniform measure are distinctly suboptimal. See \cite{adcock2023monte} for an analysis of this phenomenon.

On the other hand, Theorem \ref{thm:main-res-poly} assumes $\cV$ is a Hilbert space, whereas the lower bounds in Theorem \ref{t:lower-bounds} allow for Banach spaces. Theorem \ref{thm:main-res-poly} can be extended to Banach spaces, but with the suboptimal rate $(m/L)^{1/2(1/2-1/p)}$ \cite[Thm.\ 4.1]{adcock2023near}. It is currently unknown whether this rate can be improved.

Moreover, the computational cost of this procedure can be prohibitive. As discussed in \cite{adcock2024efficient}, combining \ef{comp-cost} with a standard estimate for $N = |\Lambda^{\mathsf{HCI}}_n|$ leads to a per-iteration costs that scales like $m^{3+\log_2(m)}$ -- i.e., subexponential, but superpolynomial in $m$. The reason for this is the need to construct the matrix $\bm{A}$ corresponding all polynomials in the large index set $\Lambda^{\mathsf{HCI}}_n$. Whether functions in $\cH(p,\mathsf{M})$ can be learned to algebraic accuracy in $m$ with polynomial computational cost in $m$ is currently an open problem. Sublinear time algorithms \cite{choi2021sparse1,choi2021sparse2} may yield a solution to this problem.

\section{DNN existence theory}\label{s:dnn-existence}

We now turn our attention to DNNs. In this section, we discuss \textit{existence theory}. Theorems of this ilk describe the \textit{expressivity} of NNs: namely, their ability to approximation functions from specific classes to a desired accuracy. We commence with a review, before showing how to establish an existence theorem for $(\bm{b},\varepsilon)$-holomorphic functions via the technique of polynomial emulation.

\subsection{Review}

Arguably, the first results on existence theory are the various \textit{universal approximation theorems} for NNs. See \cite{cybenko1989approximation,hornik1989multilayer} and, in particular, \cite{pinkus1999approximation}. These show that shallow NNs with one hidden layer can approximate any continuous function.

Unfortunately, these classical results only apply to shallow networks, and do not always give quantitative bounds on the complexity (i.e., the width) of the corresponding networks (some notable exceptions to this include \cite{mhaskar1996neural} and references therein). These issues have been investigated intensively over the last five years. Some of the first results were obtained in \cite{yarotsky2017error}. Here, explicit width and depth bounds were derived for  DNNs with the Rectified Linear Unit (ReLU) activation function for approximating Lipschitz continuous functions in the $L^{\infty}$-norm over compact sets. This inspired many subsequent works. Other activations have been studied, including hyperbolic tangents (tanh), Rectified Quadratic Units (ReQU) and Rectified Polynomial Units (RePU), and various others. And quantitative bounds have been shown for various different function classes. A partial list includes: functions in Sobolev spaces with ReLU \cite{yarotsky2017error,guhring2020error}, RePU \cite{li2020better}, tanh \cite{de-ryck2021approximation} or rational \cite{boulle2020rational} activations; piecewise smooth functions with ReLU \cite{petersen2018optimal}; $C^k$ and H\"older smooth functions with ReLU \cite{lu2021deep,schmidt-hieber2020nonparametric} or general \cite{ohn2019smooth} activations; uniformly continuous functions with ReLU activations \cite{yarotsky2018optimal}; functions in Besov spaces with ReLU \cite{opschoor2020deep,suzuki2019adaptivity}; spaces of mixed smoothness with ReLU \cite{blanchard2020representation,montanelli2019new,dung2021deep,suzuki2019adaptivity}, RePU \cite{li2020better} or smooth activations \cite{blanchard2020representation}; Gevrey functions with ReLU \cite{opschoor2020deep}, RePU \cite{opschoor2022exponential} or tanh \cite{de-ryck2021approximation} activations;  bandlimited functions with ReLU activations \cite{montanelli2021deep}; functions in Barron spaces \cite{e2021barron} with ReLU and non-ReLU activations; compositional functions \cite{poggio2017why,liang2017deep,schmidt-hieber2020nonparametric}; smooth functions on manifolds \cite{chen2022nonparametric,shaham2018provable}; finite-dimensional, analytic functions with smooth \cite{mhaskar1996neural}, ReLU \cite{e2018exponential} or RePU \cite{opschoor2022exponential} activation functions; and, most relevantly to this work, infinite-dimensional holomorphic functions with ReLU \cite{schwab2019deep,schwab2023deepa,dung2023deep}, tanh \cite{de-ryck2021approximation} or RePU \cite{schwab2023deepa} activations. See also \cite{devore2020neural,elbrachter2021deep} for reviews.

These existence theorems are based on emulating suitable classical approximation schemes with DNNs. Emulation results include (localized) Taylor polynomials \cite{de-ryck2021approximation,guhring2021approximation,liang2017deep,li2020better,lu2021deep,yarotsky2017error,schwab2019deep}, orthogonal polynomials \cite{adcock2023near,daws2019analysis,mhaskar1996neural,opschoor2022exponential,tang2019chebnet,dung2023deep,opschoor2023deep,schwab2023deepa}, rational functions \cite{boulle2020rational,telgarsky2017neural}, wavelets \cite{shaham2018provable} and general affine systems \cite{bolcskei2019optimal}, B-splines \cite{mhaskar1993approximation}, free-knot splines \cite{opschoor2020deep}, finite elements \cite{longo2023de-rham,opschoor2020deep,opschoor2023deep} and sparse grids \cite{blanchard2020representation,montanelli2019new,dung2021deep,suzuki2019adaptivity}.

To highlight some of these themes, we next describe an existence theorem for $(\bm{b},\varepsilon)$-holomorphic functions. This is achieved by first emulating Legendre polynomials with DNNs. For ease of presentation, we consider tanh DNNs. But, as we discuss, other activation functions can also be used.

\subsection{Neural network architectures}

In this and subsequent sections, we consider standard feedforward DNN architectures of the form 
\be{
\label{Phi_NN_layers}
\Phi : \bbR^n \rightarrow \bbR^k,\ \bm{z} \mapsto \Phi(\bm{z}) = \cA_{D+1} ( \sigma ( \cA_{D} ( \sigma ( \cdots \sigma ( \cA_0 (\bm{z}) ) \cdots ) ) ) ),
}
where $\cA_l : \bbR^{N_{l}} \rightarrow \bbR^{N_{l+1}}$, $l = 0,\ldots,D+1$ are affine maps and $\sigma$ is the activation function, which we assume acts componentwise. The values $\{ N_l \}^{D+1}_{l=1}$ are the widths of the hidden layers, and for convenience, we write $N_0 = n$ and $N_{D+2} = k$. Given such a DNN $\Phi$, we write 
\bes{
\mathrm{width}(\Phi) = \max \{ N_1,\ldots, N_{D+1} \},\qquad \mathrm{depth}(\Phi) = D.
}
We denote a class of DNNs of the form \eqref{Phi_NN_layers} with a fixed architecture (i.e., fixed activation function, depth and widths) as $\cN$, and define
\bes{
\mathrm{width}(\cN) = \max \{ N_1,\ldots, N_{D+1} \},\qquad \mathrm{depth}(\cN) = D.
}
Finally, since the DNNs \ef{Phi_NN_layers} take a finite input, yet the functions considered in this work take inputs in $\cU \subseteq \bbR^{\bbN}$, we also need to introduce a restriction operator. Let $\Theta \subset \bbN$ with $|\Theta| = n$. Then we define the \textit{variable restriction operator} as
\be{
\label{T-Theta_def}
\cT_{\Theta} : \bbR^{\bbN} \rightarrow \bbR^n,\ \bm{y} = (y_i)_{i \in \bbN} \rightarrow (y_i)_{i \in \Theta}.
}
If $\Theta = \{1,\ldots,n\}$ then we simply write $\cT_{\Theta} = \cT_n$.

\subsection{Emulating polynomials with DNNs: typical result}

We now present a result on emulating a finite set of multivariate Legendre polynomials with a tanh DNN.

\begin{theorem}
[Emulating multivariate Legendre polynomials]
\label{thm:dnn-emulation} 
 Let $\Lambda \subset \cF$ be a finite multi-index set and $\Theta \subset \bbN$, $|\Theta|=n$, be such that $\supp(\bm{\nu}) \subseteq \Theta$, $\forall \bm{\nu} \in \Lambda$.
Then for every $0< \delta <1$ there exists a tanh DNN $\Phi_{\Lambda,\delta}: \bbR^n \rightarrow \bbR^{|\Lambda|}$,  such that, if $\Phi_{\Lambda,\delta}(\bm{z})= (\Phi_{\bm{\nu},\delta}(\bm{z}))_{\bm{\nu} \in \Lambda}$, $\bm{z}=(z_j)_{j \in \Theta} \in \bbR^n$  and $\cT_{\Theta}$ is as in \eqref{T-Theta_def}, then  
\begin{equation*}
\|\Psi_{\bm{\nu}} -\Phi_{\bm{\nu}, \delta} \circ  \cT_{\Theta}  \|_{L^\infty(\cU)} \leq \delta, \quad \forall \bm{\nu} \in \Lambda,
\end{equation*}
where $\Psi_{\bm{\nu}}$ is the corresponding multivariate Legendre polynomial. The width and depth of this network satisfy
  \begin{align*}
  \mathrm{width}(\Phi_{\Lambda,\delta}  )  \leq c_{1}  \cdot |\Lambda|  \cdot m(\Lambda),  \quad 
  \mathrm{depth}(\Phi_{\Lambda,\delta}  ) \leq c_{2}  \cdot    \log_2(m(\Lambda)),
\end{align*}
for universal constants $c_{1},c_{2} > 0$, where $m(\Lambda) = \max_{\bm{\nu} \in \Lambda} \|\bm{\nu}\|_{1}$.
\end{theorem}

This result was shown in \cite[Thm.\ 7.4]{adcock2023near}, and is based on techniques of \cite{daws2019analysis,de-ryck2021approximation,opschoor2022exponential}. To avoid unnecessary complications, we have stated it for tanh DNNs. It also applies without changes to DNNs with the sigmoid activation, since this is obtained from tanh via shifting and scaling. The result \cite[Thm.\ 7.4]{adcock2023near} also considers the ReLU and RePU activations. As discussed therein, the proof readily extends to more general classes of activation functions. See also \cite[Rem.\ 3.9]{de-ryck2021approximation} and \cite[\S 2.4]{opschoor2022exponential} for related discussion. 

\subsection{Elements of the proof of Theorem \ref{thm:dnn-emulation}}\label{ss:emulation-proof}

The proof of Theorem \ref{thm:dnn-emulation} involves three key ingredients:
\begin{enumerate}
\item[(i)] The map $(x,y) \mapsto x y$ can be approximated by a shallow tanh NN.
\item[(ii)] Using (i), the multiplication of $d$ numbers $(x_1,\ldots,x_d) \mapsto x_1 \cdots x_d$ can be approximated by a tanh DNN of depth $c \lceil \log_2(d) \rceil$ for some constant $c>0$.
\item[(iii)] By the fundamental theorem of algebra, the multivariate Legendre polynomial $\Psi_{\bm{\nu}}$ can be expressed as a product of $\nm{\bm{\nu}}_1$ terms.
\end{enumerate}
The constructions that lead to the proofs of (i) and (ii) have become standard. In the case of tanh DNNs, they can be found in \cite[Lem.\ 3.8]{de-ryck2021approximation}. For (i), the basic idea is to use the relation $x y = ( (x+y)^2 - (x-y)^2 )/4$
in combination with a tanh DNN for approximating the map $x \mapsto x^2$. The latter is constructed via finite differences (see \cite[Lem.\ 3.1]{de-ryck2021approximation} and \cite{guhring2021approximation}).
Having shown (i), the next ingredient (ii) follows from arguments given \cite[Prop.\ 3.3]{schwab2019deep}, which formulates a DNN for multiplying $d$ numbers as a binary tree of depth $\lceil \log_2(d) \rceil$. The final ingredient (iii) follows \cite{daws2019analysis} and involves writing
\be{
\label{Psi-nu-product}
\Psi_{\bm{\nu}}(\bm{y}) = \prod_{i \in \mathrm{supp}(\bm{\nu}) } \prod^{\nu_i}_{j=1} a_{ij} (y_i - r_{ij}),\quad \forall \bm{y} \in \cU,\ \bm{\nu} \in \cF,
}
for scalars $a_{ij} , r_{ij} \in \bbR$. We deduce that $\Psi_{\bm{\nu}}(\bm{y})$ can be approximated by the composition of the affine map $\bm{y} \mapsto ( a_{ij}(y_i - r_{ij}) )_{ij}$ and the tanh DNN that approximately computes the product in \ef{Psi-nu-product}. 

The precise bounds on the architecture given in Theorem \ref{thm:dnn-emulation} are now evident. 
First, each polynomial $\Psi_{\bm{\nu}}$ involves a product of $\nm{\bm{\nu}}_1 \leq m(\Lambda)$ numbers. The depth bound now follows immediately. The width bound follows from the fact that the DNNs emulating each polynomial are stacked vertically to form the overall network $\Phi_{\Lambda,\delta}$ which simultaneously emulates all $|\Lambda|$ polynomials.

It is evident that this proof can be adapted to any other activation function for which the multiplication map in (i) can be approximated with a DNN of a quantifiable width and depth. ReLU DNNs were also considered in \cite[Thm.\ 7.4]{adcock2023near}, with ingredients (i) and (ii) being based on \cite[Prop.\ 2.6]{opschoor2022exponential}. The resulting depth and width bounds are worse than the tanh case. Conversely, a RePU network of depth one and constant width can \textit{exactly} represent the multiplication of two numbers (see, e.g., \cite[Lem.\ 2.1]{li2020better}). Therefore, a RePU network with the same width and depth bounds as in Theorem \ref{thm:dnn-emulation} can \textit{exactly} emulate the Legendre polynomials, as opposed to approximately in the tanh case. The same holds for DNNs based on \textit{rational} activation functions \cite{boulle2020rational}.

Finally, we note that the proof can be easily adapted to other polynomial systems, such as Chebyshev or, more generally, Jacobi polynomials. In the RePU case, this emulation is exact. For other activations, it is necessary to have bounds of the factors $a_{ij}$ and $r_{ij}$ in \ef{Psi-nu-product}, as the error of the approximate multiplication in (i) and (ii) depends on the size of the numbers being multiplied.

\subsection{Existence theorem for $(\bm{b},\varepsilon)$-holomorphic functions}

We now present an existence theorem for $(\bm{b},\varepsilon)$-holomorphic functions.

\thm{
[Existence theorem for DNNs]
\label{thm:existence-holo}
Let $\cV = (\bbR,\abs{\cdot})$ and $\bm{b} \in [0,\infty)^{\bbN}$ be such that $\bm{b} \in \ell^p(\bbN)$ for some $0 < p <1$. Then, for every $s \in \bbN$ there exists a class of tanh DNNs $\Phi : \bbR^s \rightarrow \bbR$ with
\bes{
\mathrm{width}(\cN) \leq c_1 \cdot s^2,\quad \mathrm{depth}(\cN) \leq c_2 \cdot \log(s), 
}
for universal constants $c_1,c_2>0$ such that the following holds. For every $f \in \cH(\bm{b})$, there exists a $\Phi \in \cN$ such that
\bes{
\nm{f - \Phi \circ \cT_{s}}_{L^2_{\varrho}(\cU)} \leq C(\bm{b},p) \cdot s^{\frac12-\frac1p}.
}
}

This result follows immediately from Theorems \ref{t:best-s-term-anchored} and \ref{thm:dnn-emulation}. The former asserts the existence of a lower set $S$ of size $|S| \leq s$ and depending on $\bm{b}$ and $p$ only which attains the desired rate of $s^{1/2-1/p}$. Then the latter asserts the existence of the DNN $\Phi_{S,\delta}$, with $\delta = s^{1/2-1/p}$, of the desired size. Note that $m(S) \leq s$ for any lower set. Indeed, if $\bm{\nu} \in S$ then so does the set $\{ \bm{\mu} \in \cF : \bm{\mu} \leq \bm{\nu} \}$. Therefore
\bes{
\nm{\bm{\nu}}_1 \leq \prod_{i \in \bbN} (\nu_i+1) = | \{ \bm{\mu} \in \cF : \bm{\mu} \leq \bm{\nu} \} | \leq | S| \leq s.
}
Theorem \ref{thm:existence-holo} only deals with scalar-valued functions. We will tackle Hilbert-valued functions in the next section when we consider learning via DNNs.

\section{Practical existence theory: near-optimal DL}\label{s:learning-proof}

While existence theory provides crucial insight into the \textit{expressivity} of DNNs, it says nothing about whether networks with similar approximation guarantees can be trained in practice and, in particular, how much training data suffices to do so.
To narrow this gap between theory and practice, in this section we develop the topic of \textit{practical existence theory}. The goal of this endeavour is to show that there exists both an architecture \textit{and} a training strategy that is similar to what is used in practice (i.e., minimizing a loss function) from which one provably learns near-optimal DNN approximations from the training data \ef{training-data-Vh}.

\subsection{Setup}

As in \S \ref{s:sparse-poly-learn}, we assume that $\cV$ is a Hilbert space with a finite-dimensional discretization $\cV_h$. Following in Remark \ref{basis-practical-soln}, we let $\{ \varphi_i \}^{K}_{i=1}$ be a basis for $\cV_h$. 
If $f : \cU \rightarrow \cV$ is the function to recover, then we may write
\bes{
f(\bm{y}) \approx  \sum^{K}_{i=1} d_{i}(\bm{y}) \varphi_i,
}
where the coefficients $d_i : \cU \rightarrow \bbR$ are scalar-valued functions. We now seek to approximate these functions using a DNN with $K$ neurons on the output layer, with the $i$th neuron corresponding to the approximation of the function $d_i$. Let $\cN$ be a class of DNNs of the form $\Phi : \bbR^n \rightarrow \bbR^K$ for some $n \in \bbN$,  
then our aim is to use the training data \ef{training-data-Vh} compute a suitable $\Phi \in \cN$ that yields an approximation $f_{\hat{\Phi}} \approx f$ defined by
\be{
\label{f-app-dnn}
f_{\Phi}(\bm{y}) : = \sum^{K}_{i=1} (\Phi \circ \cT_n(\bm{y}))_i \varphi_i,\quad \forall \bm{y} \in \cU.
}

\subsection{Practical existence theorem}

The first practical existence theorems were shown in \cite{adcock2021gap} (scalar-valued case) and \cite{adcock2021deep} (Hilbert-valued case) for holomorphic function approximation in finite dimensions with ReLU DNNs. This was extended in \cite{adcock2023near} to infinite-dimensional,  holomorphic functions with ReLU, RePU or tanh activation functions. The following result is based on \cite[Thm.\ 4.4]{adcock2023near}.

\thm{
[Practical existence theorem for DNNs]
\label{thm:main-res-dnn}
There are universal constants $c_1,c_2,c_3 \geq 1$ such that the following holds. 
Let $m \geq 3$, $0 < \epsilon < 1$ and $n = \lceil m / L \rceil$, where $L = L(m,\epsilon)$ is as in \ef{Ldef}. There there exists
\begin{enumerate}
\item a class $\cN$ of DNNs $\Phi : \bbR^n \rightarrow \bbR^K$ with tanh activation function satisfying
\bes{
\mathrm{width}(\cN) \leq c_1 \cdot m^{3 + \log_2(m)},\quad \mathrm{depth}(\cN) \leq c_2 \cdot \log(m);
}
\item a regularization function $\cJ : \cN \rightarrow [0,\infty)$ equivalent to a certain norm of the trainable parameters;
\item and a choice of regularization parameter $\lambda$ involving only $m$ and $\epsilon$,
\end{enumerate}
such that following holds for every $0 < p < 1$ and $\bm{b} \in \ell^p_{\mathsf{M}}(\bbN)$. Let $f \in \cH(\bm{b})$, $\bm{y}_1,\ldots, \bm{y}_m \sim_{\mathrm{i.i.d.}} \varrho$ and consider the training data \ef{training-data}. Then, with probability at least $1- \epsilon$, every $\gamma$-minimizer $\hat{\Phi}$, $\gamma \geq 0$, of the training problem
\begin{equation}\label{trainingprob1}
 \min_{\Phi \in \cN}\cG(\Phi),\quad \text{where } \cG(\Phi) =  \sqrt{\frac1m \sum^{m}_{i=1} \nm{f_{\Phi}(\bm{y}_i) - d_i }^2_{\cV} } +\lambda \cJ(\Phi),
\end{equation}
where $f_{\Phi}$ is as in \ef{f-app-dnn}, satisfies, with $\zeta$ as in \ef{zeta-alg-inf-def},
\be{
\label{main_err_bd}
\| f- f_{\hat{\Phi}} \|_{L^2_{\varrho}(\cU;\cV)} \leq c_3 \cdot \zeta.
}
}

In this theorem, as before, we term $\hat{\Phi}$ a \textit{$\gamma$-minimizer} of \ef{trainingprob1} if $\cG(\hat{\Phi}) \leq \min_{\Phi \in \cN} \cG(\Phi) + \gamma$.
Comparing with Theorem \ref{thm:main-res-poly}, we conclude that there is a DNN architecture and training procedure such that the resulting learned approximations achieve the same error bounds (up to possible constants) as sparse polynomial approximation based on weighted $\ell^1$-minimization. In particular, this procedure also achieves the near-optimal approximation rate $(m/L)^{\frac12-\frac1p}$.

We term Theorem \ref{thm:main-res-dnn} a \textit{practical existence theorem}. It not only asserts the existence of a DNN with a given architecture that achieves some desired rate of approximation, but demonstrates how to construct it from training data and gives generalization bounds that are explicit in the amount of training data $m$. Moreover, the training procedure is similar to standard DL procedures, in that it involves minimizing a regularized least-squares loss function. 

While Theorem \ref{thm:main-res-dnn} only considers tanh DNNs, it can be readily adapted to other activations \cite{adcock2023near}. As discussed next, all one requires to do this are variants of Theorem \ref{thm:dnn-emulation} for other activations, a topic we discussed previously in \S \ref{ss:emulation-proof}.

\subsection{The mechanism of practical existence theorems}\label{ss:pract-exist-steps}

We now give some insight into the proof of Theorem \ref{thm:main-res-dnn}, since this provides a general recipe for establishing practical existence theorems.

The overall mechanism involves not just emulating polynomials with DNNs, but also emulating the weighted $\ell^1$-minimization problem developed in \S \ref{ss:CS-setup} as a training problem.
More precisely, following \S \ref{ss:CS-setup} we proceed as follows.
\begin{enumerate}
\item[(i)] Approximate the Legendre polynomials $\{ \Psi_{\bm{\nu}} \}_{\bm{\nu} \in \Lambda}$ using DNNs.
\item[(ii)] Replace the polynomials in the matrix \ef{Ab-CS-def} the DNNs, leading to a matrix $\bm{A}' \approx \bm{A}$. Then replace $\bm{A}$ by $\bm{A}'$ in \ef{min-equiv}.
\item[(iii)] Re-cast \ef{fdef1} with $\bm{A}'$ as a training problem of the form \ef{trainingprob1}, where the unknowns $\bm{Z} \in \bbR^{N \times K}$ correspond to the trainable parameters of the DNNs.
\item[(iv)] Emulate the proof steps of Theorem \ref{thm:main-res-poly}, showing, as needed, that each step remains valid for the perturbed matrix $\bm{A}'$.
\end{enumerate}

Step (i) is accomplished straightforwardly by using the emulation result, Theorem \ref{thm:dnn-emulation}. Since $\Lambda$ is taken to be $\Lambda = \Lambda^{\mathsf{HCI}}_n$, we have $\supp(\bm{\nu}) \subseteq \{1,\ldots,n\}$ for $\bm{\nu} \in \Lambda$. Hence we set $\Theta = \{1,\ldots,n\}$ and consider a suitable parameter $\delta$ that is chosen later in the proof to balance the ensuing error terms. Let $\Phi_{\Lambda,\delta} : \bbR^n \rightarrow \bbR^{|\Lambda|}$ be the resulting DNN.

Step (ii) warrants no further discussion. Now consider Step (iii). Let $N = |\Lambda|$ as before, and define the class of DNNs
\bes{
\cN = \left \{ \Phi = \bm{Z}^{\top} \Phi_{\Lambda,\delta} : \bm{Z} \in \bbR^{N \times K} \right \}.
}
This is a class of DNNs, where only the weight matrix on the output layer is trainable. The remaining layers are nontrainable, and are handcrafted to emulate the Legendre polynomials $\{ \Psi_{\bm{\nu}} \}_{\bm{\nu} \in \Lambda}$. 

Now let $\bm{z} = (z_{\bm{\nu}_j})^{N}_{j=1} \in \cV^N_h$. Then we can associate $\bm{z}$ with its matrix of coefficients $\bm{Z} = (Z_{ij}) \in \bbR^{N \times K}_{i,j=1}$ via the relation
\be{
\label{zZrelation}
z_{\bm{\nu}_i} = \sum^{K}_{j=1} Z_{ij} \varphi_j,\quad \forall i = 1,\ldots,N,
}
and consequently with the DNN $\Phi = \bm{Z}^{\top} \Phi_{\Lambda,\delta}  \in \cN$. Using this, it is a short argument to show that
\eas{
f_{\Phi}(\bm{y}) 
 = \sum_{\bm{\nu} \in \Lambda} z_{\bm{\nu}} \Phi_{\bm{\nu},\delta} \circ \cT_{\Theta}(\bm{y}). 
}
Therefore,
\eas{
\nm{\bm{A}' \bm{z} - \bm{f}}^2_{2;\cV} &= \frac1m \sum^{m}_{i=1} \nms{ \sum_{\bm{\nu} \in \Lambda} z_{\bm{\nu}} \Phi_{\bm{\nu},\delta }(\bm{y}_i) - d_i }^2_{\cV} 
= \frac1m \sum^{m}_{i=1} \nms{f_{\Phi}(\bm{y}_i) - d_i }^2_{\cV} .
}
Now let $\cJ : \cN \rightarrow [0,\infty)$ be the regularization functional defined by 
\bes{
\cJ(\Phi) = \sum_{j=1}^{N}  u_{\bm{\nu}_j} \nms{\sum_{k=1}^K Z_{jk} \varphi_k }_{\cV} = \sum_{j=1}^{N} u_{\bm{\nu}_j}\|z_{\bm{\nu}_j} \|_{\cV}  = \nm{\bm{z}}_{1,\bm{u} ; \cV},
}
for $\Phi = \bm{Z}^{\top}  \Phi_{\Lambda,\delta} \in \cN$, where $\bm{z} = (z_{\bm{\nu}_i})^{N}_{i=1}$ is as defined by \ef{zZrelation}. We readily see that $\cJ$ is a norm over the trainable parameters. Using this and previous expression we see that \ef{fdef1} with the matrix $\bm{A}'$ in place of $\bm{A}$ can be re-cast as a training problem of the form \ef{trainingprob1}.

Finally, consider step (iv). This step is facilitated by the perturbation bound $\nm{\bm{A}-\bm{A}'}_2 \leq \sqrt{N} \delta$,
which in turn follows from a short argument via the Cauchy-Schwarz inequality:
\bes{
\nm{(\bm{A}-\bm{A}') \bm{z}}^2_2 = \frac1m \sum^{m}_{i=1} \left | \sum_{\bm{\nu} \in \Lambda} (\Psi_{\bm{\nu}}(\bm{y}_i) - \Phi_{\bm{\nu},\delta}(\bm{y}_i) ) z_{\bm{\nu}} \right |^2 \leq \left ( \sum_{\bm{\nu} \in \Lambda} \delta |z_{\bm{\nu}}| \right )^2 \leq \delta^2 N \nm{\bm{z}}^2_2.
}
The remainder of this step involves a careful modification of the proof of Theorem \ref{thm:main-res-poly} to take into account this perturbation. Theorem \ref{thm:main-res-poly} involves compressed sensing techniques, and relies on first asserting that the matrix $\bm{A}$ has a certain \textit{weighted robust Null Space Property (rNSP)} \cite{rauhut2016interpolation} (see also \cite[Chpt.\ 6]{adcock2022sparse}). This is a weighted variant of the classical rNSP, which is itself a slightly weaker condition than the better known \textit{Restricted Isometry Property (RIP)} \cite{foucart2013mathematical}. Crucially, it can be shown that the weighted rNSP is preserved under small perturbations.

We note in passing that the width and depth bounds in Theorem \ref{thm:main-res-dnn} follow quite directly from Theorem \ref{thm:dnn-emulation}. Recall that $\Lambda = \Lambda^{\mathsf{HCI}}_n$ is as in \ef{HCI-def}. A standard estimate (see \cite[Thm.\ 4.9]{kuhn2015approximation}) gives that
\be{
\label{Lambda-size}
|\Lambda| \leq \E n^{2 + \log_2(n)},\quad \forall n \in \bbN.
} 
Moreover, by definition, any $\bm{\nu} \in \Lambda$ satisfies $\nm{\bm{\nu}}_1 \leq \prod^{n-1}_{k=1} (\nu_k+1) \leq n$. Hence
\be{
\label{mLambda-bound}
m(\Lambda) \leq n.
}
The desired width and depth bounds now follow immediately from Theorem \ref{thm:dnn-emulation} and these estimates, along with the (somewhat loose) bound $n \leq m$.

\section{Epilogue}\label{s:benefits-gap}

This work has been about the approximation of smooth, infinite-dimensional functions from limited data. We close with a discussion on the benefits and consequences of practical existence theory and the gap between the \textit{handcrafted} models on which it is based and the \textit{fully trained} models used in practice.

\subsection{Scientific computing and data scarcity}

As discussed in \S \ref{ss:motivations}, parametric DE problems are often \textit{data scarce}. This is the case for many problems in scientific computing in which machine learning and, specifically, DL, is currently being applied. It also stands in stark contrast to more classical DL applications such as image classification, where datasets usually contain tens of millions of images or more. Therefore, understanding the sample complexity is crucial. Practical existence theorems show that DNNs can be efficiently learned from data.

\subsection{Potential benefits to DNNs over sparse polynomials}

Having said this, we emphasize that practical existence theorems are, at least in this work, intended primarily as \textit{theoretical} contributions. Since the strategy in Theorem \ref{thm:main-res-dnn} involves emulating the sparse polynomial approximation scheme constructed in \S \ref{s:sparse-poly-learn}, there is no benefit to implementing it over the latter.

However, in related work in inverse problems in imaging (recall the discussion in \S \ref{ss:literature}), unrolling is both used to establish practical existence theorems \textit{and} as a principled way to design DNN architectures which can then be trained as part of a DL strategy \cite{monga2021algorithm}, potentially using the theoretical weights and biases as initialization. It remains to see whether similar ideas could be effective in the parametric DE setting. Some initial work in this direction can be found in \cite{daws2019polynomial-based}.

In addition, existence theory establishes that DNNs have the capacity to approximate broad classes of functions efficiently. This is \textit{not} the case for polynomials, which fail dramatically on, for instance, discontinuous functions. This work has focused on classes of holomorphic functions, where polynomials are well suited (in fact, near optimal). DNN-based schemes have the potential to succeed on quite different function classes, something that distinguishes them from traditional methods of scientific computing. We remark in passing that discontinuous or sharp transitions arise frequently in parametric model problems, and are difficult to treat with standard methods \cite{elman2012stochastic,gunzburger2014adaptive2,gunzburger2014stochastic,jakeman2011characterization,ma2009adaptive,zhang2016hyperspherical}.

\subsection{Theorem \ref{thm:main-res-dnn} does not eliminate the theory-practice gap}

DNNs in Theorem \ref{thm:main-res-dnn} are handcrafted to emulate Legendre polynomials, with only the final layer being trained. Standard DL methods use fully trained models, where all layers are trained. Theorem \ref{thm:main-res-dnn} currently says nothing directly about this practice. However, it does lead to some insights, as we now discuss.

\subsection{Practical insights}

\subsubsection*{Width and depth bounds}

The aim of practical existence theory is to express the error in terms of the sample complexity $m$. This differs from standard existence theory in which express the error is expressed in terms of the complexity of the network, e.g., its width and depth, or its \textit{size} (number of nonzero weights and biases). Nonetheless, it is worth discussing the network complexity in Theorem \ref{thm:main-res-dnn}.

This theorem describes architectures that are much wider than they are deep. In fact, the depth grows very slowly with the number of samples $m$, like $\log(m)$. This broadly agrees with empirical insights from the application of DL in scientific computing, where it is often observed that relatively shallow networks perform well. See \cite{de-ryck2021approximation} and references therein, as well as \cite{adcock2021deep,adcock2021gap}.

However, while the depth bound in Theorem \ref{thm:main-res-dnn} is somewhat reasonable, the width bound of $\ord{m^{3+\log_2(m)}}$ grows extremely rapidly with $m$ (albeit subexponentially). 
As discussed in \S \ref{ss:pract-exist-steps}, the large estimate for the width arises, in great part, from having to emulate all polynomials in the hyperbolic cross $\Lambda^{\mathsf{HCI}}_n$, whose size, as shown in \ef{Lambda-size}, behaves like $\ord{n^{2+\log_2(n)}}$. This has nothing to do with DNN approximation itself. It stems from the sparse polynomial approximation scheme and, as discussed in \S \ref{ss:challenge1-resolve}, the need to build a finite search set $\Lambda$ outside of which the error is ensured to be small. It is unclear whether the near-optimal approximation rates of Theorem \ref{thm:main-res-poly} can be obtained with a smaller search space without further assumptions on the functions being approximated. This is an interesting prospect for future work.

\subsubsection*{Post-training pruning and sparsification}

On the other hand, it is always possible to sparsify the DNN $\hat{\Phi}$ learned in Theorem \ref{thm:main-res-dnn} after training -- a process known as \textit{pruning} (see \cite{frankle2019lottery} and references therein). This exploits the fact that the networks in $\cN$, while very wide, are sparsely connected.
This is done in much the same way as the in polynomial case (Remark \ref{rem:postprocess-sparse}). If $\hat{\bm{C}} \in \bbR^{N \times K}$ are the weights of $\hat{\Phi}$, then one first forms $\hat{\bm{c}} = (\hat{c}_{\bm{\nu}})_{\bm{\nu} \in \Lambda} \in \cV^N_h$ using \ef{zZrelation}, then computes the index set $S \subseteq \Lambda$ of the largest $n$ entries of $(\nm{\hat{c}_{\bm{\nu}}}_{\cV})_{\bm{\nu}\in \Lambda}$, and finally replaces $\hat{\Phi}$ with $\check{\Phi} = \bm{Z}^{\top}_S \Phi_{S,\delta}$, where $\bm{Z}_S \in \bbR^{n \times K}$ is formed from the rows of $\bm{Z}$ with indices in $S$. By Theorem \ref{thm:dnn-emulation} and the bound $m(S) \leq m(\Lambda) \leq n$ (recall \ef{mLambda-bound}), we have $\mathrm{width}(\check{\Phi}) \leq c_1 m^2$ and $\mathrm{depth}(\check{\Phi}) \leq c_2 \log(m)$. In particular, $\check{\Phi}$ is significantly narrow than $\hat{\Phi}$.

This suggests that methods to promote sparsity in training (see \cite{hoefler2021sparsity} and references therein) may be beneficial in practice when training fully-connected models for scientific computing problems. This requires further investigation.

\subsection{Eliminating the gap: beating the Monte Carlo rate is key}

Despite these insights, the gap between theory and practice persists. It is worth noting that standard approaches to estimating the generalization error for fully-trained models based on statistical learning theory and the \textit{bias-variance decomposition} (see, e.g., \cite{beck2022full,chen2022nonparametric,ohn2019smooth,schmidt-hieber2020nonparametric,suzuki2019adaptivity} and references therein) are not immediately applicable. These approaches use estimates for the covering number or Rademacher complexity of the relevant DNN classes. Unfortunately, they typically lead to rates that decay at best like $m^{-1/2}$. These rates are strictly slower than near-optimal rates $m^{1/2-1/p}$,  up to log factors, asserted in Theorem \ref{thm:main-res-dnn}. A major theme in parametric DEs and computational UQ is building methods that beat the \textit{Monte Carlo} rate $m^{-1/2}$ \cite[Chpt.\ 1]{adcock2022sparse}. Whether these approaches could be meaningfully combined with practical existence theorems is an interesting question for future work.

\subsection{Conclusion}

To summarize, practical existence theory is a promising way to study DNN approximation which has the potential to give new insights into the promises and challenges of data scarce applications in computational science and engineering. In addition to those outlined above, several interesting future avenues include the design of improved training methodologies and novel architectures and activation functions. Studying these areas is key for pushing the boundaries of what DNNs can achieve in scientific computing, particularly in the face of scarce data and complex computational tasks. These efforts may not only narrow the theory-to-practice gap but also unlock new DL approaches, making it more efficient, accurate, and applicable across a broader spectrum of scientific challenges and guiding the way towards more sophisticated and capable DNN models.

\bibliographystyle{abbrv}
\small
\bibliography{adcock-smooth-high-dim-refs}

\end{document}